\documentclass[12pt]{amsart}
\usepackage{amscd,amssymb,amsthm,verbatim}
\newcommand{\C}{{\mathbb C}}

\newcommand{\ch}{\operatorname{ch}}

\newcommand{\Dom}{\operatorname{Dom}}

\newcommand{\Ext}{\operatorname{Ext}}

\newcommand{\HH}{\operatorname{H}}

\newcommand{\Homeo}{\operatorname{Homeo}}

\newcommand{\Id}{\operatorname{Id}}

\newcommand{\Image}{\operatorname{Im}}
\newcommand{\Imaginary}{\operatorname{Imag}}

\newcommand{\Int}{\operatorname{int}}

\newcommand{\Isom}{\operatorname{Isom}}
\newcommand{\KK}{\operatorname{K}}
\newcommand{\KKK}{\operatorname{KK}}
\newcommand{\Ker}{\operatorname{Ker}}

\newcommand{\pt}{\operatorname{pt}}
\newcommand{\Q}{{\mathbb Q}}
\newcommand{\R}{{\mathbb R}}

\newcommand{\sign}{\operatorname{sign}}

\newcommand{\supp}{\operatorname{supp}}

\newcommand{\vol}{\operatorname{vol}}
\newcommand{\Z}{{\mathbb Z}}

\numberwithin{equation}{section}
\theoremstyle{plain}
\newtheorem{definition}[equation]{Definition}
\newtheorem{assumption}[equation]{Assumption}
\newtheorem{lemma}[equation]{Lemma}

\newtheorem{proposition}[equation]{Proposition}
\newtheorem{corollary}[equation]{Corollary}

\begin{document}

\title{Limit Sets as Examples in Noncommutative Geometry}
\author{John Lott}
\address{Department of Mathematics\\
University of Michigan\\
Ann Arbor, MI  48109-1109\\
USA} \email{lott@umich.edu}
\thanks{Research supported by NSF grant DMS-0306242}
\date{January 19, 2005}
\begin{abstract}
The fundamental group of a hyperbolic manifold acts on the limit set,
giving rise to a cross-product $C^*$-algebra.
We construct nontrivial 
K-cycles for the cross-product algebra, thereby extending
some results of Connes and Sullivan to higher dimensions.  We also
show how the Patterson-Sullivan measure on the limit set
can be interpreted as a
center-valued KMS state.
\end{abstract}
\maketitle

\section{Introduction} \label{Section 1}

If $M$ is a complete oriented $(n+1)$-dimensional
hyperbolic manifold then its fundamental
group $\Gamma$ acts on the sphere-at-infinity $S^n$ of the hyperbolic
space $H^{n+1}$. The limit set
$\Lambda$ is a closed $\Gamma$-invariant subset of $S^n$ which is the
locus for the complicated dynamics of $\Gamma$ on $S^n$. It is
self-similar and often has noninteger Hausdorff dimension.

One can associate a cross-product $C^*$-algebra $C^*(\Gamma, \Lambda)$
to the action of $\Gamma$ on $\Lambda$. It is then of interest to
see how the geometry of $M$ relates to properties of $C^*(\Gamma, \Lambda)$.
In this paper we study two aspects of this problem.  One aspect is an
interpretation of the Patterson-Sullivan measure \cite{Sullivan (1979)} in the 
framework of noncommutative geometry. The second aspect is the
construction and study of K-cycles for $C^*(\Gamma, \Lambda)$.

The Patterson-Sullivan measure is
an important tool in the study of the $\Gamma$-action on $\Lambda$.
If $x \in H^{n+1}$ then the Patterson-Sullivan measure $d\mu_x$ on $\Lambda$
describes how $\Lambda$ is seen by an
observer at $x$. 
In the first part of this paper we give an algebraic interpretation of
the Patterson-Sullivan measure.
If a $C^*$-algebra is equipped with
a one-parameter group of $*$-automorphisms then there is a notion
of a $\beta$-KMS (Kubo-Martin-Schwinger) state
on the algebra.  This notion arose from quantum statistical mechanics, 
where $\beta$ is the inverse temperature.
For each
$x \in H^{n+1}$, we construct a one-parameter group of $*$-automorphisms
of $C^*(\Gamma, \Lambda)$ and show that $d\mu_x$ gives rise to a 
$\delta(\Gamma)$-KMS state (up to normalization), 
where $\delta(\Gamma)$ is the critical
exponent of $\Gamma$. 

Putting these together for various $x$, we
obtain a picture of a field of $C^*$-algebras over $M$ with fiber
isomorphic to $C^*(\Gamma, \Lambda)$. The different copies of
$C^*(\Gamma, \Lambda)$ have different one-parameter
automorphism groups.  The global KMS state is defined on the
algebra $A$ of continuous sections of the field and takes value in 
the center $Z(A) = C(M)$. One can translate 
some geometric statements to
algebraic statements.  For example, 
if $M$ is convex-cocompact then $\delta(\Gamma)$ is the
unique $\beta$ so that there is a $\beta$-KMS state.

The bulk of the paper is concerned with constructing cycles that
represent nontrivial classes in the
K-homology $\KKK_*(C^*(\Gamma, \Lambda); \C)$ of $C^*(\Gamma, \Lambda)$,
or equivalently, in the equivariant K-homology
$\KKK_*^\Gamma(C(\Lambda); \C)$ of $C(\Lambda)$. This program was started
by Connes and Sullivan \cite[Chapter IV.3]{Connes (1994)}.
A motivation comes from the goal of doing analysis on
the self-similar set $\Lambda$. One can 
give various meanings to this phrase.  What is relevant to this paper is
the idea of Atiyah that the K-homology of a compact Hausdorff
space has cycles 
given by abstract elliptic operators on the space 
\cite{Atiyah (1970)}. This has developed into the K-homology of
$C^*$-algebras, for which we refer to the book of
Higson and Roe \cite{Higson-Roe (2000)}.
Cycles for $\KKK_*^\Gamma(C(\Lambda); \C)$ can be considered to
be something like elliptic operators on $\Lambda$.

Such cycles are pairs
$(H, F)$ satisfying certain properties, where $H$ is a Hilbert
space on which $C(\Lambda)$ and $\Gamma$ act, and $F$ is a
self-adjoint operator on $H$.
In the bounded formalism $F$ is bounded and 
commutes with the elements of $\Gamma$ up to compact operators, while 
in the unbounded formalism $F$ is generally unbounded and commutes with the
elements of $\Gamma$ up to bounded operators.

The computation of $\KKK_*^\Gamma(C(\Lambda); \C)$ can be done by
established techniques. Our goal is to find {\em explicit} and
{\em canonical} cycles
$(H, F)$ which represent nontrivial elements in 
$\KKK_*^\Gamma(C(\Lambda); \C)$. To make an analogy, a compact oriented
Riemannian manifold
has a signature class in its K-homology, but it also has
a signature {\em operator}. Clearly the study of the
signature operator leads to issues that go beyond the study of
the corresponding K-homology class.

In order to get canonical
cycles in the limit set case, we will require them
to commute with $\Gamma$ on the nose.  This is quite restrictive.
In particular, to get natural examples of
such cycles we must use the bounded formalism.
In effect, we will construct
signature-type operators on limit sets. There are two issues : first
to show that there is a nontrivial
signature-type equivariant K-homology class
on $\Lambda$, and second
to find an explicit equivariant K-cycle within the K-homology class.
Connes and Sullivan described a natural cycle
when the limit set is a quasicircle in $S^2$ and studied its
properties. As their construction
used some special features of the two-dimensional case, 
it is not immediately evident how to extend
their methods to higher dimension. 

In Section \ref{Subsection 4.1} we compute $\KKK_i^\Gamma(C(\Lambda); \C)$ in
terms of equivariant K-cohomology,
giving
$\KK^{n-i}_\Gamma(S^n, S^n - \Lambda)$. The appearance of the smooth
manifold $S^n - \Lambda$ indicates its possible relevance for constructing
K-cycles when $\Lambda \neq S^n$.

As $\Gamma$ acts conformally on $S^n$, we
construct our K-cycles in the framework of conformal geometry.
We start with the case $n = 2k$.
In Section \ref{Subsection 4.2} we consider an arbitrary 
oriented manifold $X$ of dimension $2k$, equipped with a conformal structure.
The Hilbert space $H$ 
of square-integrable $k$-forms on $X$ is conformally invariant.
We consider a certain conformally invariant
operator $F$ on $H$ that was introduced by
Connes-Sullivan-Teleman in the compact case
\cite{Connes-Sullivan-Teleman (1994)}. Under a technical assumption
(which will be satisfied in the cases of interest),
we show that $(H, F)$ gives a K-cycle for $C_0(X)$ whose K-homology
class is that of the signature operator $d + d^*$. We then prove the
invariance of the K-homology class under quasiconformal homeomorphisms
of $X$.  
This will be relevant for limit sets, as a hyperbolic manifold
has a deformation space consisting of new hyperbolic manifolds whose
dynamics on $S^n$ are conjugated to the old one by quasiconformal
homeomorphisms. 

If $\Lambda$ is the entire sphere-at-infinity $S^{2k}$ 
then the pair $(H, F)$ gives a
nontorsion class in $\KKK^\Gamma_{2k}(C(S^{2k}); \C)$.
If $\Lambda \neq S^{2k}$ then the idea will be to
sweep topological charge from $S^{2k} - \Lambda$ to $\Lambda$. 
More precisely, we have an isomorphism
$\KKK^\Gamma_{2k}(C_0(S^{2k} - \Lambda); \C) \: \cong \:
\KKK^\Gamma_{2k}(C(S^{2k}), C(\Lambda); \C)$ and a boundary map
$\KKK^\Gamma_{2k}(C(S^{2k}), C(\Lambda); \C) \rightarrow
\KKK^\Gamma_{{2k}-1}(C(\Lambda); \C)$. We can then form a
cycle in $\KKK^\Gamma_{2k-1}(C(\Lambda); \C)$ starting from the
above
cycle $(H, F)$ for $\KKK^\Gamma_{2k}(C_0(S^{2k} - \Lambda); \C)$.
Twisting the construction by $\Gamma$-equivariant
vector bundles on $S^n - \Lambda$ 
gives cycles for the
rational part of  $\KKK^\Gamma_{{2k}-1}(C(\Lambda); \C)$
represented by $\Image \left( \KK^{0}_\Gamma(S^n - \Lambda)
\rightarrow \KK^{1}_\Gamma(S^n, S^n - \Lambda) \right)$.

To make this explicit, in Section \ref{Subsection 4.3} we consider a manifold
$X$ as above equipped with a partial compactification $\overline{X}$.
Putting $\partial \overline{X} = \overline{X} - X$, for appropriate
$\overline{X}$ the pair
$(H, F)$ also gives a cycle for
$\KKK_{2k}(C_0(\overline{X}), C(\partial \overline{X}); \C)$.
The boundary map 
$\KKK_{2k}(C_0(\overline{X}), C(\partial \overline{X}); \C)
\rightarrow \KKK_{2k-1}(C(\partial \overline{X}); \C)$ was described
by Baum and Douglas in terms of Ext classes
\cite{Baum-Douglas (1991)}. 
In our case it will involve the $L^2$-harmonic $k$-forms
on $X$. If $\overline{X}$ is a smooth 
manifold-with-boundary then we show that the ensuing class
in $\Ext(C(\partial \overline{X}))$ is given by certain homomorphisms
from $C(\partial \overline{X})$ to the Calkin algebra of
a Hilbert space of
exact $k$-forms on 
$\partial \overline{X}$.
If $\overline{X}$ is the closed $2k$-ball then the Hilbert space is
the $H^{-1/2}$ Sobolev space of such forms on $S^{2k-1}$,
and is M\"obius-invariant. 

A Fuchsian group has limit set $S^{n-1} \subset S^n$. A quasiFuchsian
group is conjugate to a Fuchsian group by a
quasiconformal homeomorphism $\phi$ of $S^n$. In particular,
$\phi(S^{n-1}) = \Lambda$. In the case of a quasiFuchsian group with
$n = 2k$, we show in Section \ref{Subsection 4.4} that the element of
$\KKK^\Gamma_{{2k}-1}(C(\Lambda); \C)$ constructed by the Baum-Douglas
boundary map is
represented by the
pushforward under $\phi \big|_{S^{2k-1}}$ of the Fuchsian Ext class.
If $k = 1$ then we recover the K-homology class on a quasicircle
considered by Connes and Sullivan.
We also describe the Ext class when $M$ is
an acylindrical convex-cocompact hyperbolic $3$-manifold with incompressible
boundary, in which case $\Lambda$ is a Sierpinski curve.

Section \ref{Subsection 4.5} deals with the case when the sphere-at-infinity
$S^n$ has dimension $n = 2k-1$. If $\Lambda \neq S^{2k-1}$ then we consider
how to go from such an odd cycle on $S^{2k-1} - \Lambda$ to an even 
K-cycle on
$\Lambda$. Our discussion here is somewhat formal and uses smooth forms. 
In the case $k = 1$ we recover the K-cycle on a
Cantor set considered by Connes and Sullivan. We also describe a K-cycle
in the quasiFuchsian case and some other convex-cocompact cases.

For a quasiFuchsian limit set $\Lambda \subset S^n$, with $n$
odd or even, the K-cycle for $C(\Lambda)$ is essentially the same as
the K-cycle for $C(S^{n-1})$ in the Fuchsian case $S^{n-1} \subset S^n$, 
after pushforward by
$\phi \big|_{S^{n-1}}$.
As an example of the analytic issues concerning the K-cycle,
in Section \ref{Subsection 4.6}
we consider the subalgebra ${\mathcal A} = \phi^* C^\infty(S^n) 
\big|_{S^{n-1}}$ of $C(S^{n-1})$. We show that the Fredholm module
$({\mathcal A}, H, F)$ is $p$-summable for sufficiently large $p$.
In the case $n=2$, Connes and Sullivan showed that the infimum of such
$p$ equals $\delta(\Gamma)$. An interesting analytic question is 
how this result extends to $n > 2$.

Some related papers about limit sets are 
\cite{Anantharaman-Delaroche (1997),Anantharaman-Delaroche 
(1997b),Emerson (2003),Laca-Spielberg (1996),Spielberg (1993)}.

I thank Gilles Carron and Juha Heinonen for helpful information and
the referee for some corrections.  I thank
MSRI for its hospitality while part of this research was performed.

\section{Hyperbolic manifolds and the Patterson-Sullivan measure} 
\label{Subsection 2.1}

For background information on hyperbolic manifolds and conformal dynamics,
we refer to \cite{McMullen (1995)}. For background information
on the Patterson-Sullivan
measure, we refer to \cite{Nicholls (1989)} and
\cite{Sullivan (1979)}.

Let $\Gamma$ be a torsion-free 
discrete subgroup of $\Isom^+(H^{n+1})$, the
orientation-preserving isometries of the hyperbolic space
$H^{n+1}$. We will generally assume that $\Gamma$ is nonelementary,
i.e. not virtually abelian, although some statements will 
be clearly valid for elementary groups.
Put $M = H^{n+1}/\Gamma$, an oriented hyperbolic manifold.

We write $S^n$ for the sphere-at-infinity of $H^{n+1}$, and put
$\overline{H^{n+1}} \: = \: H^{n+1} \cup S^n$, with the topology of
the closed unit disk.
Let $\Lambda$ denote the {\em limit set} of $\Gamma$. It is the minimal 
nonempty
closed $\Gamma$-invariant subset of $S^n$.
In particular, given $x_0 \in H^{n+1}$, $\Lambda$ can be constructed as
the set of accumulation points of $x_0 \Gamma$ in $\overline{H^{n+1}}$.
The {\em domain of discontinuity} is defined to be
$\Omega \: = \: S^n \: - \: \Lambda$, an open subset of $S^n$.
There are right $\Gamma$-actions on $\Lambda$ and $\Omega$, with
the action on $\Omega$ being free and properly discontinuous.
The quotient $\Omega/\Gamma$ is called the {\em conformal boundary} of $M$.
We denote the action of $g \in \Gamma$ on $\Lambda$
by $R_g \in \Homeo(\Lambda)$. This induces a left
action of $\Gamma$ on $C(\Lambda)$, by $g \cdot f \: = \: R_g^* f$.
That is, for $g \in \Gamma$, $f \in C(\Lambda)$ 
and $\xi \in \Lambda$,
\begin{equation} \label{2.1}
(g \cdot f)(\xi) \: = \: f(\xi g).
\end{equation}

The {\em convex core} of $M$ is
the $\Gamma$-quotient of the
convex hull (in $H^{n+1}$) of $\Lambda$. The group
$\Gamma$ is {\em convex-cocompact} if the convex core
of $M$ is compact.  If $\Gamma$ is convex-cocompact then
it is Gromov-hyperbolic
and $\Lambda$ equals its Gromov boundary. 

Let $x_0$ be a basepoint in $H^{n+1}$.
The critical exponent $\delta \: = \delta(\Gamma)$ is defined by
\begin{equation} \label{2.2}
\delta \: = \: \inf \{ s \: : \: \sum_{\gamma \in \Gamma}
e^{- \: s \: d(x_0, x_0 \gamma)} \: < \: \infty \}.
\end{equation}
For each $x \in H^{n+1}$, the Patterson-Sullivan measure
$d\mu_x$ is a certain measure on $\Lambda$.
If $\Gamma$ is such that
$\sum_{\gamma \in \Gamma}
e^{- \: \delta \: d(x_0, x_0 \gamma)} \: = \: \infty$ then $d\mu_x$ is a
weak limit
\begin{equation} \label{2.3}
d\mu_x \: = \: \lim_{s \rightarrow \delta^+} \frac{\sum_{\gamma \in \Gamma}
e^{- \: s \: d(x, x_0 \gamma)} \: \delta_{x_0 \gamma}}{
\sum_{\gamma \in \Gamma}
e^{- \: s \: d(x_0, x_0 \gamma)}}
\end{equation}
of measures on $\overline{H^{n+1}}$. 
If $\sum_{\gamma \in \Gamma}
e^{- \: \delta \: d(x_0, x_0 \gamma)} \: < \: \infty$ then one
proceeds slightly differently \cite[Section 1]{Sullivan (1979)}.
% If $x_0^\prime$ is another choice of basepoint then there is a positive 
% constant
% $C_{x_0^\prime, x_0}$ so that for all $x \in H^{n+1}$, 
% $d\mu_x$ is multiplied by $C_{x_0^\prime, x_0}$.

Given $x, x^\prime \in H^{n+1}$ and $\xi \in \Lambda$, put
\begin{equation} \label{2.4}
D(x, x^\prime, \xi) \: = \: \lim_{x^{\prime \prime} \rightarrow \xi}
\left( d(x, x^{\prime \prime}) \: - \: d(x^\prime, x^{\prime \prime}) 
\right).
\end{equation}
Formally one can think of 
$D(x, x^\prime, \xi)$ as 
$d(x, \xi) \: - \: d(x^\prime, \xi)$, although the two terms do not
make individual sense.
One has
\begin{align} \label{2.5}
D(x, x^\prime, \xi) \: & = \: - \: D(x^\prime, x, \xi), \\
D(x, x^\prime, \xi) \: + \: D(x^\prime, x^{\prime \prime}, \xi) \: & = \: 
D(x, x^{\prime \prime}, \xi), \notag \\
D(x\gamma, x^\prime \gamma, \xi \gamma) \: & = \: D(x, x^\prime, \xi).
\notag 
\end{align}
One can verify from (\ref{2.3}) that
\begin{equation} \label{2.6}
d\mu_x \: = \: e^{-\delta \: D(x, x^\prime, \cdot)} \: d\mu_{x^\prime}
\end{equation}
and
\begin{equation} \label{2.7}
(R_g)_* \: d\mu_x \: = d\mu_{xg}.
\end{equation}

From (\ref{2.6}) and (\ref{2.7}),
\begin{equation} \label{2.8}
(R_g)_* \: d\mu_x \: = e^{\delta \: D(x, xg, \cdot)} \: d\mu_{x}.
\end{equation}
We note that if we have (\ref{2.7}) for a fixed $x$, and then
define $d\mu_{x^\prime}$ by (\ref{2.6}), it follows that
(\ref{2.7}) also holds for $d\mu_{x^\prime}$. We also note that
the Patterson-Sullivan measure is not a single $\Gamma$-invariant
measure.  Rather, it is a $\Gamma$-invariant conformal density in
the sense of \cite[Section 1]{Sullivan (1979)}.

\section{The cross-product $C^*$-algebra}
\label{Subsection 2.2}

The {\em algebraic cross-product} $C(\Lambda) \rtimes \Gamma$ consists of
finite formal sums $f = \sum_{g \in \Gamma} f_{g} g$, with
$f_{g} \in C(\Lambda)$.
The product of $f, f^\prime \in C(\Lambda) \rtimes \Gamma$
is given by
\begin{equation} \label{2.9}
\left( \sum_{g \in \Gamma} f_{g} g \right) 
\left( \sum_{g^\prime \in \Gamma} f^\prime_{g^\prime} g^\prime \right) \: = \:
\sum_{\gamma \in \Gamma}
\sum_{g g^\prime = \gamma} f_{g} \: (g \cdot f^\prime_{g^\prime}) \: 
\gamma,
\end{equation}
or
\begin{equation} \label{2.10}
(f f^\prime)_{\gamma}(\xi) \: = \: 
\sum_{g g^\prime = \gamma} f_{g}(\xi) \: f^\prime_{g^\prime}(\xi g).
\end{equation}
The $*$-operator is given by
\begin{equation} \label{2.11}
(f^*)_g \: = \: g \cdot \overline{f_{g^{-1}}},
\end{equation}
or
\begin{equation} \label{2.12}
(f^*)_g(\xi) \: = \: \overline{f_{g^{-1}}}(\xi g).
\end{equation}

For each $\xi \in \Lambda$, there is a $*$-homomorphism $\pi^\xi \: : \: 
C(\Lambda) \rtimes \Gamma \rightarrow B(l^2(\Gamma))$ given by
saying that for $f \: = \: \sum_{g \in \Gamma} f_g g$ and
$c \in l^2(\Gamma)$,
\begin{equation} \label{2.13}
(\pi^\xi (f) \: c)_\gamma \: = \:
\sum_{\gamma^\prime \in \Gamma} k_{\gamma, \gamma^\prime}(\xi) \: 
c_{\gamma^\prime},
\end{equation}
where
\begin{equation} \label{2.14}
k_{\gamma, \gamma^\prime}(\xi) \: = \:
f_{\gamma (\gamma^\prime)^{-1}}(\xi \gamma^{-1}).
\end{equation}
The {\em reduced cross-product $C^*$-algebra} $C^*_r(\Gamma, \Lambda)$ is the 
completion of $C(\Lambda) \rtimes \Gamma$ with respect to the norm
\begin{equation} \label{2.15}
f \rightarrow \sup_{\xi \in \Lambda} \parallel \pi^\xi \parallel_{l^2(\Gamma)}.
\end{equation}
The homomorphism $\pi^\xi$ extends to
$C^*_r(\Gamma, \Lambda)$. For $f \in C^*_r(\Gamma, \Lambda)$, 
$\pi^\xi(f)$ acts on $l^2(\Gamma)$
by a matrix
$k_{\gamma, \gamma^\prime}(\xi)$ which comes as in (\ref{2.14}) 
from a formal sum
$f \: = \: \sum_{g \in \Gamma} f_g \: g$ with each
$f_g$ in $C(\Lambda)$ 
(although if $\Gamma$ is infinite then one loses the finite support
condition when taking the completion). The product in $C^*_r(\Gamma, \Lambda)$
is given by the same formula (\ref{2.10}).

The {\em maximal cross-product $C^*$-algebra}
$C^*(\Gamma, \Lambda)$ is given by 
completing
$C(\Lambda) \rtimes \Gamma$ with respect to the supremum of
the norms of all $*$-representations
on a separable Hilbert space. 
There is an obvious homomorphism
$C^*(\Gamma, \Lambda) \rightarrow C^*_r(\Gamma, \Lambda)$.

\begin{lemma} \label{2.16}
In our case $C^*(\Gamma, \Lambda) \: = \: 
C^*_r(\Gamma, \Lambda)$. Furthermore, $C^*_r(\Gamma, \Lambda)$ is nuclear,
simple and purely infinite.
\end{lemma}
\begin{proof}
It follows from \cite[Theorem 3.1]{Spatzier-Zimmer (1991)} and
\cite[Theorem 3.37]{Anantharaman-Delaroche-Renault (2000)} that
$\Gamma$ acts topologically 
amenably on $S^n$, and hence also on
$\Lambda$. Then \cite[Proposition 6.1.8]{Anantharaman-Delaroche-Renault (2000)}
implies that $C^*(\Gamma, \Lambda) \: = \: 
C^*_r(\Gamma, \Lambda)$ and
\cite[Corollary 6.2.14]{Anantharaman-Delaroche-Renault (2000)}
implies that $C^*_r(\Gamma, \Lambda)$ is nuclear.
From \cite[Proposition 3.1]{Anantharaman-Delaroche (1997)},
$C^*_r(\Gamma, \Lambda)$ is simple and purely infinite.
(We are assuming here that $\Gamma$ is nonelementary.)
\end{proof}

Thus $C^*_r(\Gamma, \Lambda)$ is a Kirchberg algebra
\cite{Rordam (2002)}. In addition, it lies in the so-called
bootstrap class ${\mathcal N}$, as follows for example from 
\cite[Section 10]{Tu (1997)}. Thus
$C^*_r(\Gamma, \Lambda)$
falls into a class of $C^*$-algebras that can be classified by
their K-theory.

\section{An automorphism group and 
a positive functional on $C^*_r(\Gamma, \Lambda)$}
\label{Subsection 3.1}

In this section, for each
$x \in H^{n+1}$, we construct a corresponding
one-parameter group of $*$-automorphisms
of $C^*(\Gamma, \Lambda)$. We show that the Patterson-Sullivan measure
$d\mu_x$ gives rise to a 
$\delta(\Gamma)$-KMS state (up to normalization).

Propositions \ref{3.2} and \ref{3.7} of the present
section are special cases of general results about
quasi-invariant measures and KMS states
\cite[Chapter II.5]{Renault (1980)}.  We include the proofs, which are quite 
direct in our case, for completeness.

Fix $x \in H^{n+1}$.
Given $t \in \R$ and $f \in C^*_r(\Gamma, \Lambda)$, put
\begin{equation} \label{3.1}
(\alpha_t f)_{g}(\xi) \: = \:
e^{itD(x, xg^{-1}, \xi)} \: f_{g}(\xi).
\end{equation}

\begin{proposition} \label{3.2}
$\alpha$ is a strongly-continuous one-parameter group of $*$-automorphisms of 
$C^*_r(\Gamma, \Lambda)$.
\end{proposition}
\begin{proof}
For $f \in C(\Lambda) \rtimes \Gamma$, 
the kernel $k^t$ corresponding to $\alpha_t f$ is
\begin{align} \label{3.2.5}
k^t_{\gamma, \gamma^\prime}(\xi) \: & = \:
(\alpha_t f)_{\gamma (\gamma^\prime)^{-1}}(\xi \gamma^{-1})
\: = \:
e^{itD(x, x\gamma^\prime \gamma^{-1}, \xi \gamma^{-1})} \: 
f_{\gamma (\gamma^\prime)^{-1}}(\xi \gamma^{-1}) \: = \: 
e^{itD(x, x\gamma^\prime \gamma^{-1}, \xi \gamma^{-1})} \: 
k_{\gamma, \gamma^\prime}(\xi) \\
& = \:
e^{itD(x \gamma, x\gamma^\prime, \xi)} \: 
k_{\gamma, \gamma^\prime}(\xi) \: = \:
e^{it(D(x \gamma, x, \xi) - D(x\gamma^\prime, x, \xi))} 
\: k_{\gamma, \gamma^\prime}(\xi). \notag
\end{align}
Thus $\pi^{\xi}(\alpha_t f) \: = \: U(t, \xi) \: \pi^{\xi}(f) \: U(t, \xi)^{-1}$, where
$U(t, \xi)$ is the unitary operator that acts on $c \in l^2(\Gamma)$ by
\begin{equation} \label{3.3}
(U(t, \xi) c)_g \: = \: e^{it D(xg, x, \xi)} \: c_g.
\end{equation} 
This shows that if $f \in C^*_r(\Gamma, \Lambda)$ then
$\alpha_t f \in C^*_r(\Gamma, \Lambda)$, and that
$\alpha_t f$ is strongly-continuous in $t$.

Given $f, f^\prime \in  C^*_r(\Gamma, \Lambda)$,
\begin{align} \label{3.4}
(\alpha_t (f f^\prime))_{g}(\xi) \: & = \:
e^{itD(x, xg^{-1}, \xi)} \: (f f^\prime)_{g}(\xi) \: = \:
e^{itD(x, xg^{-1}, \xi)} \: 
\sum_{\gamma \gamma^\prime = g} f_{\gamma}(\xi) \: 
f^\prime_{\gamma^\prime}(\xi \gamma) \\
& = \:
\sum_{\gamma \gamma^\prime = g} e^{itD(x, x\gamma^{-1}, \xi)} \:
f_{\gamma}(\xi) \: e^{itD(x \gamma^{-1}, x g^{-1}, \xi)}
f^\prime_{\gamma^\prime}(\xi \gamma) \notag \\
& = \:
\sum_{\gamma \gamma^\prime = g} e^{itD(x, x\gamma^{-1}, \xi)} \:
f_{\gamma}(\xi) \: e^{itD(x, x (\gamma^\prime)^{-1}, \xi \gamma)}
f^\prime_{\gamma^\prime}(\xi \gamma) \: = \:
((\alpha_t f) (\alpha_t f^\prime))_{g}(\xi). \notag
\end{align}
Thus $\alpha_t (f f^\prime) \: = \: (\alpha_t (f)) (\alpha_t (f^\prime))$.
Next, given $f \in C^*_r(\Gamma, \Lambda)$,
\begin{align} \label{3.5}
(\alpha_t f)^*_g(\xi) \: & = \: \overline{
e^{itD(x,xg,\xi g)} \: f_{g^{-1}}(\xi g)} \: = \:
e^{-itD(x,xg,\xi g)} \: \overline{f_{g^{-1}}(\xi g)} \: = \:
e^{-itD(xg^{-1},x,\xi)} \: \overline{f_{g^{-1}}(\xi g)} \\
& = \: e^{itD(x, xg^{-1},\xi)} \: \overline{f_{g^{-1}}(\xi g)} \: = \:
(\alpha_t f^*)_g(\xi). \notag
\end{align}
Thus $(\alpha_t (f))^* \: = \: \alpha_t (f^*)$.
This shows that $\alpha_t$ is a $*$-automorphism of $C^*_r(\Gamma, \Lambda)$.
Finally, it is clear that for $t, t^\prime \in \R$,
$\alpha_t \circ \alpha_{t^\prime} \: = \: \alpha_{t + t^\prime}$.
\end{proof}

Define a positive functional 
$\tau \: : \: C^*_r(\Gamma, \Lambda) \rightarrow \C$ by
\begin{equation} \label{3.6}
\tau (f) \: = \:
\int_{\Lambda} f_{e} \: d\mu_x. 
\end{equation}
It may not be a state, as $d\mu_x$ may not be a probability measure.
(See Lemma \ref{3.25}. One could imagine normalizing 
$d\mu_x$ by dividing it by its mass, but this would cause further
complications.)

For background on KMS states, we refer to
\cite[Chapter 8.12]{Pedersen (1979)}. 
We now show that $\tau$ satisfies the KMS condition.

\begin{proposition} \label{3.7}
Given $f, f^\prime \in C^*_r(\Gamma, \Lambda)$, and $t \in \R$, put
\begin{equation} \label{3.8}
F(t) \: = \: \tau(f \: \alpha_t(f^\prime)). 
\end{equation}
Then $F$ has a continuous bounded continuation to $\{z \in \C \: : \:
0 \le \Imaginary(z) \le \delta\}$ that is analytic in
$\{z \in \C \: : \:
0 < \Imaginary(z) < \delta\}$, with
\begin{equation} \label{3.9}
F(t+i\delta) \: = \: \tau(\alpha_t(f^\prime) \: f). 
\end{equation}
\end{proposition}
\begin{proof}
From \cite[Proposition 8.12.3]{Pedersen (1979)}, it is enough to show
that (\ref{3.9}) holds when $f^\prime \in C(\Lambda) \rtimes \Gamma$.
In this case,
\begin{align} \label{3.10}
F(t) \: & = \: \int_{\Lambda} (f (\alpha_tf^\prime))_{e}(\xi) \: d\mu_x(\xi) \:
= \: \int_{\Lambda} \sum_{g \in \Gamma} f_g(\xi) \: 
(\alpha_t f^\prime)_{g^{-1}}(\xi g) \: d\mu_x(\xi) \\ 
& = \:
\int_{\Lambda} \sum_{g \in \Gamma} f_g(\xi) \: 
e^{it D(x, xg, \xi g)} \:
f^\prime_{g^{-1}}(\xi g) \: d\mu_x(\xi) \notag \\
& = \:
\int_{\Lambda} \sum_{g \in \Gamma} f_g(\xi) \: 
e^{it D(xg^{-1}, x, \xi)} \:
f^\prime_{g^{-1}}(\xi g) \: d\mu_x(\xi). \notag
\end{align}
Then
\begin{align} \label{3.11}
F(t+i\delta) \: 
& = \:
\int_{\Lambda} \sum_{g \in \Gamma} f_g(\xi) \: 
e^{it D(xg^{-1},x, \xi)} \:
f^\prime_{g^{-1}}(\xi g) \: e^{-\delta D(xg^{-1},x, \xi)} \:
d\mu_x(\xi) \\
& = \:
\int_{\Lambda} \sum_{g \in \Gamma} f_g(\xi) \: 
e^{it D(xg^{-1},x, \xi)} \:
f^\prime_{g^{-1}}(\xi g) \: 
d\mu_{xg^{-1}}(\xi) \notag \\
& = \:
\int_{\Lambda} \sum_{g \in \Gamma} f_g(\xi) \: 
e^{it D(xg^{-1}, x,\xi)} \:
f^\prime_{g^{-1}}(\xi g) \: 
d\mu_{x}(\xi g) \notag \\
& = \:
\int_{\Lambda} \sum_{g \in \Gamma} f_g(\xi g^{-1}) \: 
e^{it D(xg^{-1},x, \xi g^{-1})} \:
f^\prime_{g^{-1}}(\xi) \: 
d\mu_{x}(\xi) \notag \\
& = \:
\int_{\Lambda} \sum_{g \in \Gamma}e^{it D(xg,x, \xi g)} \:
f^\prime_{g}(\xi) \: f_{g^{-1}}(\xi g) \: 
d\mu_{x}(\xi) \notag \\
& = \:
\int_{\Lambda} \sum_{g \in \Gamma}e^{it D(x,xg^{-1}, \xi)} \:
f^\prime_{g}(\xi) \: f_{g^{-1}}(\xi g) \: 
d\mu_{x}(\xi) \notag \\
& = \: \int_{\Lambda} \sum_{g \in \Gamma} (\alpha_t f^\prime)_g(\xi) \: 
f_{g^{-1}}(\xi g) \: d\mu_x(\xi) \notag \\ 
& = \: \int_{\Lambda} ( (\alpha_tf^\prime) f)_{e}(\xi) \: d\mu_x(\xi) 
 \: = \: \tau(\alpha_t(f^\prime) \: f). \notag
\end{align}
This proves the claim.
\end{proof}

\section{Center-valued KMS state} \label{Subsection 3.2}

In this section we allow the point $x \in H^{n+1}$ to vary.
We construct a field of $C^*$-algebras over $M$, each 
isomorphic to $C^*(\Gamma, \Lambda)$. 
The global KMS state is defined on the
algebra $A$ of continuous sections of the field and takes value in 
the center $Z(A) = C(M)$. We translate some statements about the
conformal dynamics of $\Gamma$ on $S^n$ to
algebraic statements about the KMS state on $A$.

Let $C(H^{n+1}, C^*_r(\Gamma, \Lambda))$ denote the continuous 
maps from $H^{n+1}$ to $C^*_r(\Gamma, \Lambda)$. 
We write an element of $C(H^{n+1}, C^*_r(\Gamma, \Lambda))$
as $F \: \equiv
\: \sum_{g \in \Gamma} F_{x, g} g$, with $F_{x,g} \in C(\Lambda)$.
Then $\Gamma$ acts by automorphisms on 
$C(H^{n+1}, C^*_r(\Gamma, \Lambda))$, by the formula
\begin{equation} \label{3.12}
(\gamma \cdot F)_{x,g} \: = \: R_\gamma^* F_{x \gamma, \gamma^{-1} g \gamma}.
\end{equation}

Define
an $1$-parameter group of $*$-automorphisms
${\mathcal A}_t$ of
$C(H^{n+1}, C^*_r(\Gamma, \Lambda))$ by
\begin{equation} \label{3.13}
({\mathcal A}_t F)_{x, g}(\xi) \: = \:  
e^{itD(x, xg^{-1}, \xi)} \: F_{x,g}(\xi).
\end{equation}
\begin{lemma} \label{3.14}
${\mathcal A}_t$ is $\Gamma$-equivariant.
\end{lemma}
\begin{proof}
Given $\gamma \in \Gamma$ and 
$F \in C(H^{n+1}, C^*_r(\Gamma, \Lambda))$,
\begin{align} \label{3.15}
({\mathcal A}_t (\gamma \cdot F))_{x, g}(\xi) \: & = \:  
e^{itD(x, xg^{-1}, \xi)} \: F_{x \gamma, \gamma^{-1} g \gamma}(\xi \gamma)
\: = \: 
e^{itD(x\gamma, x \gamma \gamma^{-1} g^{-1} \gamma, \xi \gamma)} 
\: F_{x \gamma, \gamma^{-1} g \gamma}(\xi \gamma) \\
& = \:
(\gamma \cdot ({\mathcal A}_t F))_{x,g} (\xi). \notag
\end{align}
\end{proof}

We write the positive functional $\tau$ of (\ref{3.6}) as $\tau_x$. For
$F \in C(H^{n+1}, C^*_r(\Gamma, \Lambda))$, define 
${\mathcal T} (F) \in C(H^{n+1})$ by
\begin{equation} \label{3.16}
({\mathcal T} (F))(x) \: = \:  
\tau_x \left(\sum_{g \in \Gamma} F_{x, g} g \right) \: = \:
\int_\Lambda F_{x,e}(\xi) \: d\mu_x(\xi). 
\end{equation}

\begin{lemma} \label{3.17}
${\mathcal T}$ is $\Gamma$-equivariant.
\end{lemma}
\begin{proof}
Given $\gamma \in \Gamma$ and 
$F \in  C(H^{n+1}, C^*_r(\Gamma, \Lambda))$,
\begin{align} \label{3.18}
(R_\gamma^* ({\mathcal T} (F)))(x) \: & = \: 
({\mathcal T} (F))(x\gamma) \: = \:  
\tau_{x\gamma} \left(\sum_{g \in \Gamma} F_{x\gamma, g} g \right)
\: = \: \int_\Lambda  F_{x\gamma, e} \: d\mu_{x\gamma} \\
& = \: \int_\Lambda  F_{x\gamma, e} \: (R_\gamma)_* d\mu_{x}
\: = \: \int_\Lambda (R_\gamma)^* F_{x\gamma, e} \: d\mu_{x}
\: = \: \int_\Lambda (\gamma \cdot F)_{x, e} \: d\mu_{x} \notag \\
& = \: ({\mathcal T} (\gamma \cdot F))(x). \notag
\end{align}
\end{proof}

Let ${A}$ be the $\Gamma$-invariant subspace 
$\left(  C(H^{n+1}, C^*_r(\Gamma, \Lambda)) \right)^{\Gamma}$.
Then ${A}$ consists of the continuous sections of a field of
$C^*$-algebras over $M$ in the sense of
\cite[Definition 10.3.1]{Dixmier (1977)}, with each fiber $A_m$ 
isomorphic to
$ C^*_r(\Gamma, \Lambda)$. The center of $A$ is $Z(A) \: = \:
C(M)$.
By Lemma \ref{3.14}, the automorphisms
${\mathcal A}_t$ restrict to a $1$-parameter group ${\mathcal B}_t$ of
$*$-automorphisms of $A$.

By Lemma \ref{3.17}, the map 
${\mathcal T}$ restricts to a map
${\mathcal S} \: : \: A \rightarrow Z(A)$.
For $F, F^\prime \in A$, put
\begin{equation} \label{3.19}
{\mathcal F}(t) \: = \: {\mathcal S}(F \: {\mathcal B}_t(F^\prime)). 
\end{equation}
As in Proposition \ref{3.7}, ${\mathcal F}$ has a continuous extension
to $\{z \in \C \: : \:
0 \le \Imaginary(z) \le \delta\}$ that is analytic
(in the sense of \cite[Definition 3.30]{Rudin (1973)}) in
$\{z \in \C \: : \:
0 < \Imaginary(z) < \delta\}$, with
\begin{equation} \label{3.20}
{\mathcal F}(t+i\delta) \: = \: {\mathcal S}({\mathcal B}_t(F^\prime) \: F). 
\end{equation}

\begin{lemma} \label{3.21}
For all $\sigma \in Z(A)$ and $F \in A$,
${\mathcal S}(\sigma F) \: = \: \sigma {\mathcal S}(F)$. 
\end{lemma}

We will call a linear map ${\mathcal S} \: : \: A \rightarrow Z(A)$ satisfying
the preceding properties a center-valued $\delta$-KMS state for the
pair $(A, {\mathcal B}_t)$, or just a
$\delta$-KMS state.
We do not require that ${\mathcal S}(1_A)$ be $1$.

\begin{proposition} \label{3.22}
a. If $\Gamma$ is convex-cocompact then the pair 
$(A, {\mathcal B}_t)$ has a $\delta(\Gamma)$-KMS state, and this is the only
$\beta$ for which $(A, {\mathcal B}_t)$ has a $\beta$-KMS state.
Furthermore, the KMS-state is unique up to multiplication by positive
elements of $Z(A)$.\\
b. If $\Gamma$ is not convex-cocompact then for each 
$\beta \in [\delta(\Gamma), \infty)$ 
the pair 
$(A, {\mathcal B}_t)$ has a $\beta$-KMS state.\\
c. If $\Gamma$ is not convex-cocompact and has no parabolic elements
then the set of
$\beta$ for which $(A, {\mathcal B}_t)$ has a $\beta$-KMS state
is $[\delta(\Gamma), \infty)$.
\end{proposition}
\begin{proof}
Existence of $\beta$ : For all $\Gamma$, the Patterson-Sullivan measure gives
rise to a $\delta(\Gamma)$-KMS state on $(A, {\mathcal B}_t)$.
Fix $x \in H^{n+1}$.
From  \cite[Theorem 2.19(i)]{Sullivan (1987)}, if $\Gamma$ is not
convex-cocompact then for each
$\beta \in [\delta(\Gamma), \infty)$, 
there is a positive measure $d\nu_x$ on $\Lambda$ satisfying
\begin{equation} \label{3.23}
(R_g)_* \: d\nu_x \: = e^{\beta \: D(x, xg, \cdot)} \: d\nu_{x}.
\end{equation}
Given such a measure, for $x^\prime \in H^{n+1}$, define 
$d\nu_{x^\prime}$ by
\begin{equation} \label{3.24}
d\nu_{x^\prime} \: = \: e^{\beta \: D(x, x^\prime, \cdot)} \: d\nu_{x}.
\end{equation}
Then we can form a $\beta$-KMS state for the pair $(A, {\mathcal B}_t)$
in the same way as with the Patterson-Sullivan measure. \\
Uniqueness of $\beta$ : Suppose that $\Gamma$ has no parabolic elements.
Fix $x \in H^{n+1}$. Consider the cross-product groupoid $G = \Lambda
\rtimes \Gamma$.
Define the cocycle $c(\xi, g) \: = \: D(x, xg, \xi)$. Suppose that
$\xi g \: = \: g$ and
$c(\xi, g) \: = \: 0$. Take an upper half-plane model for $H^{n+1}$ in which
$\xi$ is the point at infinity. Then the hyperbolic element $g$ 
translates by a signed length $d(g)$ in the $(n+1)$-th coordinate
(along with a possible rotation in the other coordinates), and 
$|D(x, xg, \xi)| \: = \: |d(g)|$. It follows that $g$ is the identity
element of $\Gamma$. Thus the subgroupoid $c^{-1}(0)$ is principal.
 
Suppose that we have a $\beta$-KMS state for the pair 
$(A, {\mathcal B}_t)$. From \cite[Proposition 3.2]{Kumjian-Renault (2003)},
the KMS state arises from a positive measure $d\nu_x$ on $\Lambda$
which satisfies (\ref{3.23}).
Then from \cite[Theorem 2.19]{Sullivan (1987)},
if $\Gamma$ is convex-cocompact then $\beta \: = \: \delta(\Gamma)$, while if
$\Gamma$ is not cocompact then $\beta \in [\delta(\Gamma), \infty)$.
Furthermore, if $\Gamma$ is convex-cocompact then $d\nu_x$ is 
proportionate to the Patterson-Sullivan measure $d\mu_x$.
\end{proof}

\begin{lemma} \label{3.25}
${\mathcal S}(1)$ is a positive eigenfunction of $\triangle_M$ with
eigenvalue $\delta(\Gamma) \: (n-\delta(\Gamma))$.
\end{lemma}
\begin{proof}
The function $\Phi$ on $H^{n+1}$, given by setting
$\Phi(x)$ to be the mass of $d\mu_x$, is the pullback to $H^{n+1}$ of
a positive eigenfunction $\phi$ of $\triangle_M$ with eigenvalue
$\delta(\Gamma) \: (n - \delta(\Gamma))$
\cite[Theorem 28]{Sullivan (1979)}.
\end{proof}

In general, ${\mathcal S}(1)$ is not bounded on $M$.

\begin{lemma} \label{3.26}
If $\Gamma$ is convex-cocompact then
${\mathcal S}(1) \in C_0(M)$.
\end{lemma}
\begin{proof}
With reference to the proof of Lemma \ref{3.25}, if $\Gamma$ is 
convex-cocompact then
\cite[Theorem 2.13(a)]{Sullivan (1987)} implies that $\phi \in C_0(M)$,
from which the result follows.
\end{proof}

In the rest of this section we assume that $\Gamma$ is convex-cocompact.
Let $\pi_m \: : \: A \rightarrow A_m$ be the homomorphism from $A$ to the
fiber over $m \in M$.
Let $A_0$ be the subalgebra of $A$ consisting of
elements $a$ so that the function $m \rightarrow 
\parallel \pi_m(a) \parallel$ lies in $C_0(M)$.
Then $A_0$ is the $C^*$-algebra associated to the continuous field of
$C^*$-algebras on $M$, in the sense of 
\cite[Section 10.4.1]{Dixmier (1977)}. From Lemma \ref{3.26}, the map
${\mathcal S} \: : \: A \rightarrow Z(A)$ restricts to a map
${\mathcal S}_0 \: : \: A_0 \rightarrow Z(A_0)$, for which
(\ref{3.19}) and (\ref{3.20}) again hold.  Also, 
for all $\sigma \in Z(A_0)$ and $F \in A_0$,
${\mathcal S}_0(\sigma F) \: = \: \sigma {\mathcal S}_0(F)$.

\section{K-homology of the cross-product algebra} \label{Subsection 4.1}

In this section we compute $\KKK_i^\Gamma(C(\Lambda); \C)$ in
terms of the equivariant K-cohomology, in the sense of the
Borel construction, of the pair $(S^n, \Omega)$.

We let $\KK^*(\cdot, \cdot)$ denote the representable (i.e. homotopy-invariant)
K-cohomology of a 
topological pair 
\cite[Chapter 7.68, Remark in Chapter 8.43, Chapter 11]{Switzer (1975)}.
We let $\KK_*(\cdot)$ denote
the unreduced Steenrod K-homology of a compact metric space
\cite[p. 161]{Ferry (1995)},
\cite{Kahn-Kaminker-Schochet (1977)}.
Put $\overline{M} \: = \: (H^{n+1} \cup \Omega)/\Gamma$, so
$\overline{M} \: = \: M \cup \partial \overline{M}$, where
$\partial \overline{M}$ is the conformal boundary.

For background on analytic K-homology and (equivariant) KK-theory,
we refer to \cite{Higson-Roe (2000)} and
\cite{Blackadar (1998)}.
We recall that $\KKK^\Gamma_i(C(\Lambda); \C)$ is isomorphic to
$\KKK_i(C^*(\Gamma, \Lambda); \C)$ \cite[Theorem 20.2.7]{Blackadar (1998)}. 
We wish to
compute $\KKK_i^\Gamma(C(\Lambda); \C)$ in term of classical
homotopy-invariant topology (as opposed to
proper homotopy invariance). 

If $X$ and $A \subset X$ are manifolds then
the relative K-group $\KK^0(X, A)$ has generators given by
virtual vector bundles on $X$ that are trivialized over $A$, and
similarly for $\KK^1(X, A)$. We let
$\KK^*_\Gamma (X, A)$ denote the relative K-theory of the
Borel construction, e.g. $\KK^*_\Gamma (S^n, \Omega) \: = \:
\KK^*((E\Gamma \times S^n)/\Gamma, (E\Gamma \times \Omega)/\Gamma)$.
A model for $E\Gamma$ is $H^{n+1}$. There is a 
$\Gamma$-equivariant diffeomorphism $SH^{n+1} \rightarrow H^{n+1} \times
S^n$ that sends a unit vector $\widehat{v}$ at a point $x \in H^{n+1}$ to the pair
$(x, \xi)$, where $\xi$ is the point on the sphere-at-infinity hit by the
geodesic starting at $x$ with initial vector $\widehat{v}$. Passing to
$\Gamma$-quotients gives a diffeomorphism
$SM \rightarrow (E\Gamma \times S^n)/\Gamma$.
The subspace $(E\Gamma \times \Omega)/\Gamma$ can be identified
with the unit tangent vectors $v \in SM$ with the property that
the geodesic generated by $v$ goes out the conformal boundary
$\partial \overline{M}$. We note that $(E\Gamma \times \Omega)/\Gamma$
is homotopy-equivalent to $\partial \overline{M}$.

\begin{proposition} \label{4.1}
$\KKK_i(C(\Lambda); \C) \: \cong \: \KK^{n-i}(S^n,
\Omega)$ and
$\KKK^\Gamma_i(C(\Lambda); \C) \: \cong \: \KK^{n-i}_\Gamma(S^n,
\Omega)$.
\end{proposition}
\begin{proof}
We have $\KKK_i(C(\Lambda); \C) \cong \KK_{i}(\Lambda)$
\cite[Theorem C]{Kahn-Kaminker-Schochet (1977)}.
By Alexander duality \cite[Theorem B]{Kahn-Kaminker-Schochet (1977)}, 
\begin{equation} \label{4.2}
\KK_{i}(\Lambda) \: \cong \:
\KK^{n-i}(S^n, \Omega).
\end{equation}
(The statement of \cite[Theorem B]{Kahn-Kaminker-Schochet (1977)}
is in terms of reduced homology and cohomology, but is equivalent
to (\ref{4.2}) if $\Omega$ is nonempty. The case when $\Omega$ is
empty is more standard
\cite[Theorem 14.11]{Switzer (1975)}.)

There is a spectral sequence to compute
$\KKK^\Gamma_{-i}(C(\Lambda); \C)$, with differential of degree $+1$ and
$E_2$-term given by
$E_2^{p,q} \: = \: \HH^p(\Gamma, \KKK_{-q}(C(\Lambda); \C))
\: = \:\HH^p(\Gamma, \KK_{-q}(\Lambda))$
\cite[Theorem 6]{Kasparov (1984)},
\cite[p. 199]{Kasparov (1988)}. As $B\Gamma$ has a model that is a
finite-dimensional CW-complex, there is no problem with convergence of
the spectral sequence.
By (\ref{4.2}), $\KK_{-q}(\Lambda) \: \cong \:
\KK^{n+q}(S^n, \Omega)$.
Then
$E_2^{p,q} \: \cong \:
\HH^p(\Gamma, \KK^{n+q}(S^n, \Omega))$. 
This will be the same as $E_2$-term of the Leray spectral
sequence \cite[Theorem 15.27, Remarks 2 and
3 on p. 351-352]{Switzer (1975)}
to compute $\KK^{n+i}_\Gamma(S^n, \Omega)$ from the
fibration $((E \Gamma \times S^n)/\Gamma, (E \Gamma \times \Omega)/\Gamma)
\rightarrow B\Gamma$, with the same
differentials. Changing the sign of $i$ gives
the claim.
\end{proof}

The significance of Proposition \ref{4.1} is that when $\Omega \neq 
\emptyset$, 
it indicates that it should be possible to construct elements of 
$\KKK_*^\Gamma( C(\Lambda); \C)$ by means of the smooth manifold
$\Omega$.
More precisely, we have an isomorphism
$\KKK^\Gamma_{n}(C_0(\Omega); \C) \: \cong \:
\KKK^\Gamma_{n}(C(S^{n}), C(\Lambda); \C)$ and a boundary map
$\KKK^\Gamma_{n}(C(S^{n}), C(\Lambda); \C) \rightarrow
\KKK^\Gamma_{{n}-1}(C(\Lambda); \C)$. We can then
start with an explicit cycle $(H, F)$ for 
$\KKK^\Gamma_{n}(C_0(\Omega); \C)$ and follow these maps to
construct the corresponding
cycle in $\KKK^\Gamma_{{n}-1}(C(\Lambda); \C)$.

If $\Lambda \: = \: S^n$ then the signature
class $\sigma \in \KKK_n(C(S^n); \C)$ satisfies
$\sigma \: = \: C_n \: [S^n]$, where
$ [S^n] \in \KKK_n(C(S^n); \C)$ is the fundamental K-homology
class, represented by the
Dirac operator, and $C_n$ is a power of $2$. 
Under the isomorphism (\ref{4.2}), $\sigma$ goes over
to $* \sigma \: = \: C_n \:  [1] \in \KK^0(S^n)$.
Applying the Chern character gives
$\ch(* \sigma) \: = \: C_n \cdot 1 \in \HH^0(S^n; \Q)$.

There is a natural transformation $f \: : \:
\KK^*(X, A) \otimes \Q \rightarrow
\KK^*(X, A; \Q)$, where the right-hand-side is K-theory with
coefficients. For general topological spaces,
 $f$ need not be injective or surjective.
If $X$ and $A$ are finite-dimensional $CW$-complexes
then the Atiyah-Hirzebruch
spectral sequence implies that $f$ is injective and
has dense image in the sense that the annihilator of
$\Image(f)$, in the dual space $\left( \KK^*(X,A; \Q) \right)^*$, vanishes.
(Note that tensoring with $\Q$ does not commute with arbitrary direct 
products.)
If in addition 
$\KK^*(X, A; \Q)$ is finite-dimensional then $f$ is an isomorphism.
From the proof of Proposition \ref{4.1}
there is an injective map 
$\KKK^\Gamma_i(C(\Lambda); \C) \otimes \Q \rightarrow
\KK^{n-i}_\Gamma(S^n,
\Omega; \Q)$ with dense image,
which is an isomorphism when the right-hand-side is
finite-dimensional.

The Chern character gives an isomorphism between
$\KK^*(X, A; \Q)$ and
$\HH^*(X, A; \Q)$, after $2$-periodization of the latter, and
similarly for $\KK^*_\Gamma(X, A; \Q)$. One can compute 
$\HH^*_\Gamma(S^n, \Omega; \Q)$ using the Leray spectral sequence,
with $E_2$-term $E_2^{p,q} \: = \: \HH^p(\Gamma; \HH^q(S^n, \Omega; \Q))$.
If $\Lambda \: = \: S^n$ then 
$E_2^{0,0} \: = \: \HH^0(\Gamma; \HH^0(S^n; \Q)) \: = \:
\HH^0(\Gamma; \Q) \: = \: \Q$. This term is unaffected by 
the differentials of the
spectral sequence, and so it passes to the limit.
In particular, the element
$C_n \cdot 1 \in \HH^0(S^n; \Q)$ is $\Gamma$-invariant and gives a 
nonzero element of $\HH^0_\Gamma(S^n; \Q) \: = \: \Q$. Hence there is
a corresponding element of 
$\KK^0_\Gamma(S^n; \Q)$.

If $\Lambda \neq S^n$ and $n > 1$ then the exact sequence
\begin{equation}
0 \rightarrow \HH^0(S^n, \Omega; \Q) \rightarrow
\HH^0(S^n; \Q) \rightarrow
\HH^0(\Omega; \Q) \rightarrow
\HH^1(S^n, \Omega; \Q) \rightarrow
\HH^1(S^n; \Q) \rightarrow \ldots
\end{equation}
implies that
$\HH^0(S^n, \Omega; \Q) \: = \: 0$ and 
$\HH^1(S^n, \Omega; \Q) \: = \: \Q^{|\pi_0(\Omega)|}/\Q$.
Then the $E_2^{p,0}$-term of the spectral sequence
for $\HH^*_\Gamma(S^n, \Omega; \Q)$ vanishes, and
the $E_2^{0,1}$-term is
$\HH^0(\Gamma; \HH^1(S^n; \Q)) \: = \:
\HH^0(\Gamma; \Q^{|\pi_0(\Omega)|}/\Q) \: \cong \:
\Q^{|\pi_0(\partial \overline{M})|}/\Q$. 
This term is unaffected by
the differentials of the
spectral sequence, and so it passes to the limit
to give a contribution to $\HH^1_\Gamma(S^n, \Omega; \Q)$. There is a
corresponding component of $\KK^1_\Gamma(S^n, \Omega; \Q)$.

\section{An even K-cycle on a manifold} \label{Subsection 4.2}

In this section we consider an arbitrary 
oriented manifold $X$ of dimension $2k$, equipped with a conformal structure.
The Hilbert space $H$ 
of square-integrable $k$-forms on $X$ is conformally invariant.
We consider a certain conformally invariant
operator $F$ that was introduced by
Connes-Sullivan-Teleman in the compact case
\cite{Connes-Sullivan-Teleman (1994)}. Under a technical assumption,
we show that $(H, \gamma, F)$ gives a K-cycle for $C_0(X)$ whose K-homology
class is that of the signature operator $d + d^*$. We then show the
invariance of the K-homology class under quasiconformal homeomorphisms. 

As a short digression, let us discuss why we use the operator $F$. 
It is well-known that the bounded K-cycle 
$\left( L^2(X; \Lambda^*), \frac{d \: + \: d^*}{\sqrt{1 + \triangle}} \right)$
represents a nontrivial
class in $K_{2k}(C_0(X))$.
In the case $X \: = \: S^{2k}$ , equipped with the 
action of a discrete group $\Gamma$ by M\"obius transformations, this operator 
gives rise to an element of $\KK_{2k}^\Gamma(C(S^{2k}); \C)$ 
\cite{Kasparov (1984)}, but at the price of making some modifications.
Namely,
there is a natural action of $\Gamma$ on the
$L^2$-forms on $S^{2k}$ which is unitary on 
$L^2(S^{2k}; \Lambda^k)$ but is nonunitary on $L^2(S^{2k}; \Lambda^*)$
(as we are using a Riemannian structure). One has
to modify the $\Gamma$-action in order to make it unitary.  After doing so,
the $\Gamma$-action commutes with 
$\frac{d \: + \: d^*}{\sqrt{1 + \triangle}}$ {\em up to compact operators}.
In later sections we will take
$X \: = \: \Omega \: = \: S^{2k} - \Lambda$, on which the relevant group
$\Gamma$ 
acts conformally.  We want a K-cycle that {\em commutes}
with $\Gamma$. The Connes-Sullivan-Teleman operator is well-suited for
this purpose. In addition,
the conformal {\em invariance} of the Connes-Sullivan-Teleman
operator will lead to the quasiconformal invariance of its
K-homology class. This will be important when we consider quasiconformal
deformations of $\Gamma$-actions.

For notation, if $X$ is a Riemannian manifold then
we let $L^2(X; \Lambda^q)$ denote
the square-integrable $q$-forms on $X$, and similarly for
$L^p(X; \Lambda^q)$, $L^p_c(X; \Lambda^q)$, $C^\infty(X; \Lambda^q)$, 
$C^\infty_c(X; \Lambda^q)$ and $H^s(X; \Lambda^q)$, where the $c$-subscript
denotes compact support.

\subsection{Some conformally-invariant constructions}
\label{Subsubsection 4.2.1}

In this subsection we define the operator $F$ and introduce the technical
Assumption \ref{4.13}.

As for the role of Assumption \ref{4.13},
if $X$ is compact then one can use a
pseudodifferential calculus to see that 
$(H, \gamma, F)$ gives a K-cycle for $C(X)$. If $X$ is
noncompact then there is a local pseudodifferential calculus
on $X$, but it will be
insufficient to verify the K-cycle conditions. Instead we use
finite-propagation-speed arguments for Dirac-type
operators.  Assumption \ref{4.13} effectively arises in
interpolating between our operator $F$ and the Dirac-type operator
$D \: = \: d + d^*$.
 
Let $X$ be an oriented $2k$-dimensional manifold with a given conformal 
class $[g]$ of Riemannian metrics.
\begin{lemma} \label{4.4}
There is a complete Riemannian metric in the conformal class.
\end{lemma}
\begin{proof}
Without loss of generality, we may assume that $X$ is connected.
Choose a Riemannian metric $g_0$ in the conformal class.
There is an exhaustion $K_0 \subset K_1 \subset \ldots$ of $X$ by
smooth compact manifolds-with-boundary, with $K_i \subset
\Int(K_{i+1})$.  For $i > 1$, choose a nonnegative smooth function
$\phi_i$ with $\supp(\phi_i) \subset \Int(K_{i+1}) \: - \: K_{i-2}$ so that
for any path $\{\gamma_i(t)\}_{t \in [0,1]}$ 
from $\partial K_{i-1}$ to $\partial
K_i$, $\int_0^1 e^{\phi_i(\gamma_i(t))} 
\: g_0(\gamma_i^\prime, \gamma_i^\prime)
^{1/2} \: dt \: \ge 1$. Put
$\phi \: = \: \sum_i \phi_i$. Then $g \: = \: e^{2\phi} \: g_0$ is complete.
\end{proof}

We now make some constructions that are independent of the choice of
the complete Riemannian metric $g$ in the conformal class $[g]$.
Consider the 
complex Hilbert space $H \: = \: L^2(X; \Lambda^k)$ 
of square-integrable
$k$-forms on $X$, with its conformally-invariant inner product. 
There is an obvious action of $C_0(X)$ on $H$.
Let $\gamma$ be the conformally-invariant 
$\Z_2$-grading operator on $H$ given by
\begin{equation} \label{4.5}
\gamma \: = \: i^{k} \: *.
\end{equation}
Let $H \: = \: H_+ \oplus H_-$ be the corresponding orthogonal
decomposition. There are operators
\begin{equation} \label{4.6}
d \: : \: C^\infty_c(X; \Lambda^{k-1}) \rightarrow
C^\infty_c(X; \Lambda^{k})
\end{equation}
and
\begin{equation} \label{4.7}
d^* \: : \: C^\infty_c(X; \Lambda^{k+1}) \rightarrow
C^\infty_c(X; \Lambda^{k}).
\end{equation}
Then 
\begin{equation}
\Image(d^*) \: = \: \gamma \: \Image(d).
\end{equation}
There is a conformally-invariant orthogonal decomposition
\begin{equation} \label{4.8}
H \: = \: \overline{\Image(d)} \: \oplus \: \overline{\Image(d^*)} \: \oplus
{\mathcal H},
\end{equation}
where
\begin{equation} \label{4.9}
{\mathcal H} \: = \: \{ \omega \in H \cap C^\infty(X; \Lambda^k) \: : \:
d\omega \: = \: d^* \omega \: = \: 0\}.
\end{equation}
Furthermore, ${\mathcal H}$ is an orthogonal direct sum 
${\mathcal H}_+ \oplus {\mathcal H}_-$ of its self-dual and
anti-self-dual subspaces. 

We note that the normed vector space
$L_c^{\frac{2k}{k-1}}(X; \Lambda^{k-1})$ is conformally-invariant.

\begin{lemma} \label{4.10}
$\overline{\Image(d)}$ equals the closure of the image of $d$ on
$\{ \eta \in L_c^{\frac{2k}{k-1}}(X; \Lambda^{k-1}) \: : \: d\eta \in  
L^{2}(X; \Lambda^{k})\}$.
\end{lemma}
\begin{proof}
Clearly $\overline{\Image(d)}$ is contained in the closure of the image of 
$d$ on
$\{ \eta \in L_c^{\frac{2k}{k-1}}(X; \Lambda^{k-1}) \: : \: d\eta \in  
L^{2}(X; \Lambda^{k})\}$. Conversely, suppose that $\eta
\in L_c^{\frac{2k}{k-1}}(X; \Lambda^{k-1})$ has $d\eta \in  
L^{2}(X; \Lambda^{k})$. Let $\rho \in C^\infty_c(\R)$ be an even
function with support in $[-1, 1]$ and
$\int_\R \rho(s) \: ds \: = \: 1$. 
Put $\triangle \: = \: d d^* \: + \: d^* d$. 
For $\epsilon \: > \: 0$, put
\begin{equation} \label{4.11}
\widehat{\rho}(\epsilon^2 \triangle) \: = \:
\int_\R e^{is\epsilon(d+d^*)} \: \rho(s) \: ds \: = \:
\int_\R \cos(s\epsilon \sqrt{\triangle}) \: \rho(s) \: ds.
\end{equation}
By elliptic theory, 
$\widehat{\rho}(\epsilon^2 \triangle) \eta \in C^\infty(X; \Lambda^{k-1})$.
By finite propagation speed arguments
\cite[Proposition 10.3.1]{Higson-Roe (2000)}, the support of 
$\widehat{\rho}(\epsilon^2 \triangle) \eta$ lies within distance $\epsilon$ of
the essential support of $\eta$, so 
$\widehat{\rho}(\epsilon^2 \triangle) \eta \in 
C^\infty_c(X; \Lambda^{k-1})$. Finally, by the functional calculus,
$\lim_{\epsilon \rightarrow 0} d \left(
\widehat{\rho}(\epsilon^2 \triangle) \eta \right) \: = \:
\lim_{\epsilon \rightarrow 0}  
\widehat{\rho}(\epsilon^2 \triangle) d \eta \: =
d \eta$ in $L^{2}(X; \Lambda^{k})$.
\end{proof}

Define $F \in B(H)$ by
\begin{equation} \label{4.12}
F(\omega) \: = \: 
\begin{cases}
\omega & \text{ if } \omega \in \overline{\Image(d)}, \\
- \: \omega & \text{ if } \omega \in \overline{\Image(d^*)}, \\
0 & \text{ if } \omega \in {\mathcal H}.
\end{cases}
\end{equation}
Then $F^* \: = \: F$ and $F$ anticommutes with $\gamma$.

\begin{assumption} \label{4.13}
There is a complete Riemannian metric in the conformal class 
such that for each $\omega \in \overline{\Image(d)}$, there is
an $\eta \in L^2(X; \Lambda^{k-1})$ with $d\eta \: = \: \omega$.
\end{assumption}

We do not know if Assumption \ref{4.13} is really necessary for what
follows, but it is required for our proofs.  It is equivalent to
saying that there is a gap away from zero in the spectrum of the
Laplacian on $L^2(X; \Omega^k)$
\cite[Proposition 1.2]{Lott (1996)}. \\ \\
% An alternative assumption that does not make reference to the Riemannian
% metric would be to say that for each 
% $\omega \in \overline{\Image(d)}$, there is
% an $\eta \in L^{\frac{2k}{k-1}}(X; \Lambda^{k-1})$ 
% with $d\eta \: = \: \omega$. We do not know if there is a relationship
%between this alternative assumption and Assumption \ref{4.13}.
{\bf Example 1 : } Assumption \ref{4.13} is
satisfied for the conformal class of the unit ball in $\R^{2k}$, by 
taking the hyperbolic metric. More generally, it is satisfied when $X$ is
the interior of a compact manifold-with-boundary $\overline{X}$, and
the conformal class comes from a smooth Riemannian metric $g_0$ on 
$\overline{X}$. One can
see this by using the complete asymptotically-hyperbolic
metric on $X$ given by
$g \: = \: \rho^{-2} \: g_0$, where near the boundary 
$\partial{\overline{X}}$, $\rho \in C^\infty(\overline{X})$ 
equals the distance function to
the boundary with respect to $g_0$. Then the essential spectrum of
the $k$-form Laplacian on $X$ will be the same as that of the
essential spectrum of the $k$-form Laplacian on $H^{2k}$, which
has a gap away from zero.
 \\ \\
{\bf Example 2 : } Assumption \ref{4.13} is 
satisfied for the conformal class
of the standard Euclidean metric on $\R^{2k}$. Consider a radially
symmetric metric on
$\R^{2k}$ of the form 
$g \: = \: \sigma^2(r) \: (dr^2 \: + \: r^2 \: d\theta^2)$,
where $\sigma \: : \: (0, \infty) \rightarrow (0, \infty)$ is a smooth
function satisfying
\begin{equation}
\sigma(r) \: = \: 
\begin{cases}
1 & \text{ if } r < 1, \\
\frac{1}{r \ln{r}} & \text{ if } r > 2.
\end{cases}
\end{equation}
From \cite[Theorem 2.2]{Donnelly-Xavier (1984)}, the essential
spectrum of the $k$-form Laplacian on $(\R^{2k}, g)$ is bounded
below by a positive constant. (In the case $k = 1$,
$(\R^2, g)$ has a hyperbolic cusp at infinity.) \\ \\
{\bf Example 3 : } Suppose that a discrete group $\Gamma$ acts properly and
cocompactly on $X$. Considering metrics on $X$ that pullback from the
orbifold $X/\Gamma$, whether or not Assumption \ref{4.13} is satisfied 
for these metrics is
topological, i.e. independent of the metric on $X/\Gamma$.

\subsection{A conformally-invariant K-cycle}
\label{Subsubsection 4.2.2}

In this subsection, under Assumption \ref{4.13}, we show that 
$(H, \gamma, F)$ gives
a K-cycle for $C_0(X)$ whose K-homology class is that of the signature
operator $d + d^*$.

For notation,
if $H$ is a Hilbert space then we denote the bounded operators on $H$ by
$B(H)$, the compact operators on $H$ by $K(H)$ and the Calkin algebra by
$Q(H) \: = \: B(H)/K(H)$.
We recall that a cycle for $\KKK_0(C_0(X); \C)$
is given by a triple $(H, \gamma, F)$ where \\
1. $H$ is a separable Hilbert space with $\Z_2$-grading operator $\gamma \in
B(H)$,\\
2. There is a $*$-homomorphism $C_0(X) \rightarrow B(H)$ and\\
3. $F \in B(H)$ is such that $F \gamma \: + \: \gamma F \: = \: 0$ and
for all $a \in C_0(X)$, we have
$a (F^2 - I) \in K(H)$, $a (F - F^*) \in K(H)$ and
$[F, a] \in K(H)$. \\

We now consider the triple $(H, \gamma, F)$ of Section 
\ref{Subsubsection 4.2.1}.
We let $P_{\overline{\Image(d)}}$, $P_{\overline{\Image(d^*)}}$ and
$P_{\mathcal H}$ denote orthogonal projections onto
$\overline{\Image(d)}$, $\overline{\Image(d^*)}$ and
${\mathcal H}$, respectively. We let $G$ denote the Green's operator for
$\triangle$ on $L^2(X; \Lambda^k)$, so $\triangle G \: = \: G \triangle 
\: = \: I \: - \: P_{\mathcal H}$.

\begin{proposition} \label{4.14}
For all $a \in C_0(X)$, $a(F^2 - I)$ is compact.
\end{proposition}
\begin{proof}
We may assume that $a \in C^\infty_c(X)$.  This
is because for any $a \in C_0(X)$, there is a sequence
$\{a_i\}_{i=1}^\infty$ in $C^\infty_c(X)$ with
$\lim_{i \rightarrow \infty} a_i \: = \: a$ in the sup norm. Then
$a(F^2 - I)$ will be the norm limit of the compact operators
$a_i (F^2 - I)$, and hence compact.

We have
$I - F^2 \: = \: P_{\mathcal H}$. Let $K$ be the support of $a$.
Choose a complete Riemannian metric $g$ in the given conformal class.
Applying G\r{a}rding's inequality \cite[10.4.4]{Higson-Roe (2000)}
with $D \: = \: d \: + \: d^*$,
there is a $c > 0$ so that for all $\omega \in H$,
\begin{equation} \label{4.15}
c \parallel P_{\mathcal H} \omega \parallel_{H^1(K; \Lambda^k)} \: \le
\: \parallel P_{\mathcal H} \omega \parallel_{L^2(M; \Lambda^k)} \: \le
\: \parallel \omega \parallel_{L^2(M; \Lambda^k)}.
\end{equation}
It follows that the map $\omega \rightarrow a (P_{\mathcal H} \omega) \big|_K$
is bounded from $L^2(M; \Lambda^k)$ to $H^1(K; \Lambda^k)$. By Rellich's
Lemma \cite[10.4.3]{Higson-Roe (2000)}, the inclusion map from
$H^1(K; \Lambda^k)$ to $L^2(M; \Lambda^k)$ is compact. The proposition
follows.
\end{proof}

\begin{proposition} \label{4.16}
If Assumption \ref{4.13} is satisfied then
for all $a \in C_0(X)$, $[F, a]$ is compact. 
\end{proposition}
\begin{proof}
It is enough to prove the proposition for $a \in C^\infty_c(X)$. 
We may assume that $a$ is real.
Write the action of $a$ on $H$ as a $(3 \times 3)$-matrix with
respect to the decomposition (\ref{4.8}). Then we must show that its
off-diagonal entries are compact.  By the self-adjointness of $a$,
it is enough to show that
$(I \: - \: P_{\overline{\Image(d)}}) \: a P_{\overline{\Image(d)}} \: : \: 
\overline{\Image(d)} \rightarrow 
\overline{\Image(d^*)} \: \oplus
{\mathcal H}$
and 
$(I \: - \: P_{\overline{\Image(d^*)}}) \: a P_{\overline{\Image(d^*)}} 
\: : \: 
\overline{\Image(d^*)} \rightarrow 
\overline{\Image(d)} \: \oplus
{\mathcal H}$
are compact. 

Given $\eta \in C^\infty_c(X; \Lambda^{k-1})$,
\begin{align} \label{4.17}
a \: d\eta \: & = \: d(a\eta) \: - \: da \wedge \eta \\
& = \: d(a\eta) \: - \: da \wedge (P_{\mathcal H} \eta \: + \:
d G d^* \eta \: + \: d^* G d \eta) \notag \\
& = \: d(a(\eta \: - \: P_{\mathcal H} \eta \: - \:
d G d^* \eta)) \: - \: da \wedge G^{1/2} d^* G^{1/2} d \eta. \notag 
\end{align}
Thus
\begin{equation} \label{4.18}
(I \: - \: P_{\overline{\Image(d)}}) \: a P_{\overline{\Image(d)}}
\: = \: - \: 
(I \: - \: P_{\overline{\Image(d)}}) \: da \wedge G^{1/2} d^* G^{1/2} 
P_{\overline{\Image(d)}}.
\end{equation}
As $d^* G^{1/2}$ is bounded, 
to show that $(I \: - \: P_{\overline{\Image(d)}}) \: a
P_{\overline{\Image(d)}}$ is compact, it suffices to show that
$da \wedge G^{1/2} \: : \: (\overline{\Image(d^*)} \subset
L^2(X; \Lambda^{k-1})) \rightarrow
L^2(X; \Lambda^k)$ is compact. Put $D \: = \: d \: + \: d^*$, so
$D^2 \: = \: \triangle$.
By Assumption \ref{4.13}, there is an even
function $\rho \in C_0(\R)$ so that when acting on
$\overline{\Image(d^*)} \subset
L^2(X; \Lambda^{k-1})$, we have
$G^{1/2} \: = \: \rho(D)$. We can assume that
$\rho(x) \: = \: \frac{1}{|x|}$ for $|x|$
large. The compactness now follows from the
fact that $da \wedge \rho(D) \: : \: L^2(X; \Lambda^{k-1}) \rightarrow
L^2(X; \Lambda^{k})$ is compact
\cite[Proposition 10.5.2]{Higson-Roe (2000)}.

Let $(da)_{\sharp}$ denote the vector field that is dual to $da$, with
respect to $g$.
Given $\eta \in C^\infty_c(X; \Lambda^{k+1})$,
\begin{align} \label{4.19}
a \: d^* \eta \: & = \: d^*(a\eta) \: + \: i_{(da)_\sharp} \eta \\
& = \: d^*(a\eta) \: + \: i_{(da)_\sharp} (P_{\mathcal H} \eta \: + \:
d G d^* \eta \: + \: d^* G d \eta) \notag \\
& = \: d^*(a(\eta \: - \: P_{\mathcal H} \eta \: - \:
d^* G d \eta)) \: + \: i_{(da)_\sharp} G^{1/2} d G^{1/2} d^* \eta. \notag 
\end{align}
Thus
\begin{equation} \label{4.20}
(I \: - \: P_{\overline{\Image(d^*)}}) \: a P_{\overline{\Image(d^*)}}
\: = \: 
(I \: - \: P_{\overline{\Image(d^*)}}) \: i_{(da)_\sharp} G^{1/2} d G^{1/2} 
P_{\overline{\Image(d^*)}}.
\end{equation}
Following the previous line of proof, we conclude that
$(I \: - \: P_{\overline{\Image(d^*)}}) \: a P_{\overline{\Image(d^*)}} 
\: : \: 
\overline{\Image(d^*)} \rightarrow 
\overline{\Image(d)} \: \oplus
{\mathcal H}$
is compact. 
\end{proof}

Thus the triple $(H, \gamma, F)$ is a cycle for
$\KKK_{0}(C_0(X); \C) \cong \KKK_{2k}(C_0(X); \C)$.
 We extend $\gamma$ to the usual $\Z_2$-grading
on $L^2(X; \Lambda^*)$.

\begin{proposition} \label{4.21}
If Assumption \ref{4.13} is satisfied then
the cycles $(H, \gamma, F)$ and
$(L^2(X; \Lambda^*), \gamma, \frac{d+d^*}{\sqrt{1+\triangle}})$ 
represent
the same class in
$\KKK_{2k}(C_0(X); \C)$.
\end{proposition}
\begin{proof}
Define $\overline{F} \in B(L^2(X; \Lambda^*))$ by
\begin{equation} \label{4.22}
\overline{F} \omega \: = \:
\begin{cases}
\omega & \text{ if } \omega \in L^2(X; \Omega^j), \: j \: < \: k, \\
F\omega & \text{ if } \omega \in L^2(X; \Omega^k),\\
- \: \omega & \text{ if } \omega \in L^2(X; \Omega^j), \: j \: > \: k.
\end{cases}
\end{equation}
Then $\overline{F}$ anticommutes with $\gamma$, and the cycle
$(L^2(X; \Lambda^*), \gamma, \overline{F})$ differs from
$(H, \gamma, F)$ by the addition of a degenerate cycle.  Hence they
define the same class in $\KKK_{2k}(C_0(X); \C)$. Now 
$\overline{F}$ commutes with $\frac{d+d^*}{\sqrt{1+\triangle}}$, so
it
anticommutes with $i \: \gamma \: \frac{d+d^*}{\sqrt{1+\triangle}}$. Then the
cycles with $F_t \: = \: \cos(t) \: \overline{F} \: + \: i \: \sin(t)
\: \gamma \frac{d+d^*}{\sqrt{1+\triangle}}$, $t \in [0, \frac{\pi}{2}]$,
homotop from $(L^2(X; \Lambda^*), \gamma, \overline{F})$ to
$(L^2(X; \Lambda^*), \gamma, i \: \gamma \: \frac{d+d^*}{\sqrt{1+\triangle}})$.
Finally, the cycles with $F_t \: = \: (i \: \gamma \: \cos(t) \: + \: \sin(t))
\: \frac{d+d^*}{\sqrt{1+\triangle}}$, $t \in [0, \frac{\pi}{2}]$, homotop 
from $(L^2(X; \Lambda^*), \gamma, i \: \gamma \:
\frac{d+d^*}{\sqrt{1+\triangle}})$ to 
$(L^2(X; \Lambda^*), \gamma, \frac{d+d^*}{\sqrt{1+\triangle}})$.
The proposition follows.
\end{proof}
\noindent
{\bf Remark : } If $X$ is compact then Proposition \ref{4.21} 
was previously proved
in \cite[p. 677]{Connes-Sullivan-Teleman (1994)} by a different argument.

\subsection{Quasiconformal invariance}
\label{Subsubsection 4.2.3}

In this subsection we show that the K-homology class of $(H, \gamma, F)$ 
is invariant
under quasiconformal homeomorphisms of $X$.

\begin{proposition} \label{4.23}
If $\phi : X_1 \rightarrow X_2$ is an orientation-preserving
$K$-quasiconformal homeomorphism, for some $K < \infty$,
and $X_1$ and $X_2$ satisfy Assumption \ref{4.13}, then
$\phi_* [(H_1, \gamma_1, F_1)] \: = \:  [(H_2, \gamma_2, F_2)]$ in
$\KKK_{2k}(C_0(X_2); \C)$.
\end{proposition}
\begin{proof}
The pushforward $\phi_* [(H_1, \gamma_1, F_1)] \in
\KKK_{2k}(C_0(X_2); \C)$ is represented by a K-cycle using 
$H_1$, $\gamma_1$ and $F_1$,
where
$C_0(X_2)$ acts on $H_1$ via the pullback $\phi^* \: : \:
C_0(X_2) \rightarrow C_0(X_1)$. As $\phi$ is
$K$-quasiconformal, $(\phi^{-1})^* H_1$ and $H_2$ are the
same as topological
vector spaces. By naturality, we can 
represent $\phi_* [(H_1, \gamma_1, F_1)]$ by letting $C_0(X_2)$ act on 
$(\phi^{-1})^* H_1$, equipped with the transported operator
$(\phi^{-1})^* F_1$. 
From Lemma \ref{4.10}, $(\phi^{-1})^* \overline{\Image(d)} \: = \:
\overline{\Image(d)}$.
Then
$(\phi^{-1})^* F_1$ is the operator constructed using $d$ and the
transported grading operator
$(\phi^{-1})^* \gamma_1$. Hence it suffices to work on a fixed manifold $X$
and consider two conformal structures that are $K$-quasiconformal.
Equivalently, we can
consider the corresponding grading operators $\gamma_1$ and
$\gamma_2$ \cite[Lemma 2.3]{Donaldson-Sullivan (1989)}.

There is a measurable bundle homomorphism $\mu_+ \: : \:
\Lambda^k_- \rightarrow \Lambda^k_+$ with 
$\sup_{x \in X} |\mu_+(x)| \: < \: 1$ so that
if $\mu \: = \:
\begin{pmatrix}
0 & \mu_+ \\
\mu_+^* & 0
\end{pmatrix}$
then
$\gamma_2 \: = \: (1+\mu) \gamma_1 (1+\mu)^{-1}$
\cite[Section 4$\alpha$]{Connes-Sullivan-Teleman (1994)},
\cite[Section 2(i)]{Donaldson-Sullivan (1989)}.
For $t \in [0, 1]$, put $\gamma(t) \: = \: (1+t\mu) \gamma_1 (1+t\mu)^{-1}$.
The corresponding inner product space has
\begin{equation} \label{4.24}
\langle \omega_1, \omega_2 \rangle (t) \: = \:
\langle \omega_1, (1-t\mu) \: (1+t\mu)^{-1} \: \omega_2 \rangle (0).
\end{equation}
The operator $F(t)$ is
one on $\overline{\Image(d)}$, minus one on 
$\gamma(t) \overline{\Image(d)}$ and zero on
$\left(\overline{\Image(d)} \oplus \gamma(t) \overline{\Image(d)} 
\right)^\perp$.

The Hilbert spaces $\{H(t)\}_{t \in [0,1]}$ form a Hilbert
$C([0, 1])$-module. They all have the same underlying topological vector
space. We claim that the operators $\{F(t)\}_{t \in [0,1]}$
are norm-continuous in $t$. For this, it suffices to show that
the projection operators $P_{\overline{\Image(d)}}$ and
$P_{\overline{\Image(d^*)}}$ are norm-continuous in $t$. 
As $\overline{\Image(d)}$ is independent of $t$,
\cite[Lemma 6.2]{Hilsum (1985)} implies that
$P_{\overline{\Image(d)}}$ is norm-continuous in $t$.
As $\Ker(d) \: = \: \overline{\Image(d)}
\oplus {\mathcal H}$ is independent of $t$, it also follows from 
\cite[Lemma 6.2]{Hilsum (1985)} that $P_{\overline{\Image(d)}} \: + \:
P_{\mathcal H}$ is norm-continuous in $t$. Then
$P_{\overline{\Image(d^*)}} \: = \: I \: - \: 
P_{\overline{\Image(d)}} \: - \:
P_{\mathcal H}$ is norm-continuous in $t$. 

The operators 
$\gamma(t)$ are also norm-continuous in $t$.
In order to show that $\{(H(t), \gamma(t), F(t))\}_{t \in [0, 1]}$
is a homotopy of $K$-cycles, it now suffices to show that 
for all $a \in C_0(X)$, $[F(t), a]$ and
$a (F(t)^2 - 1)$ are compact operators.
We may assume that
$a \in C^\infty_c(X)$. From Propositions \ref{4.14} and \ref{4.16},
$[F(0), a]$ and
$a (F(0)^2 - 1)$ are compact. Using the fact that
$\frac{d}{dt} d^* \: = \: \left[
\frac{d\gamma}{dt} \: \gamma^{-1}, d^* \right]$, one can compute that
\begin{align} \label{4.25}
\frac{d}{dt} P_{\overline{\Image(d)}} \: & = \: - \:
 P_{\overline{\Image(d)}} \: \frac{d\gamma}{dt} \: \gamma^{-1} \:
(I \: - \: P_{\overline{\Image(d)}}), \\
\frac{d}{dt} P_{\overline{\Image(d^*)}} \: & = \:
(I \: - \: P_{\overline{\Image(d^*)}}) \: 
\frac{d\gamma}{dt} \: \gamma^{-1} \: P_{\overline{\Image(d)^*}}, \notag \\
\frac{d}{dt} P_{\mathcal H} \: & = \: - \:
P_{\mathcal H} \: 
\frac{d\gamma}{dt} \: \gamma^{-1} \: P_{\overline{\Image(d^*)}}
\: + \: P_{\overline{\Image(d)}}  \: 
\frac{d\gamma}{dt} \: \gamma^{-1} \: P_{\mathcal H}. \notag
\end{align}

To compute $\frac{d}{dt} [F(t), a]$, it suffices to compute
$\frac{d}{dt} [P_{\overline{\Image(d)}}, a]$ and 
$\frac{d}{dt} [P_{\overline{\Image(d^*)}}, a]$.
Now
\begin{align} \label{4.26}
\frac{d}{dt} \left[ P_{\overline{\Image(d)}}, a \right] \: & = \: - \:
\left[
 P_{\overline{\Image(d)}} \: \frac{d\gamma}{dt} \: \gamma^{-1} \:
(I \: - \: P_{\overline{\Image(d)}}), a \right] \\
& = \:
\: - \:
\left[
 P_{\overline{\Image(d)}}, a \right] \: \frac{d\gamma}{dt} \: \gamma^{-1} \:
(I \: - \: P_{\overline{\Image(d)}}) \: - \:
 P_{\overline{\Image(d)}} \: \frac{d\gamma}{dt} \: \gamma^{-1} \:
\left[
(I \: - \: P_{\overline{\Image(d)}}), a \right] \notag \\
& = \:
\: - \:
\left[
 P_{\overline{\Image(d)}}, a \right] \: \frac{d\gamma}{dt} \: \gamma^{-1} \:
(I \: - \: P_{\overline{\Image(d)}}) \: + \:
 P_{\overline{\Image(d)}} \: \frac{d\gamma}{dt} \: \gamma^{-1} \:
\left[ P_{\overline{\Image(d)}}, a \right]. \notag
\end{align}
From the proof of Proposition \ref{4.16}, 
at $t = 0$, $\left[ P_{\overline{\Image(d)}}(0), a \right]$ is
compact. 
From (\ref{4.26}), we can write $\left[ P_{\overline{\Image(d)}}(t), a \right]
\: = \: U(t) \: \left[ P_{\overline{\Image(d)}}(0), a \right] \: V(t)$,
where $U(0) \: = \: V(0) \: = \: I$ and
\begin{align} \label{4.27}
\frac{dU}{dt} \: & = \: P_{\overline{\Image(d)}}(t) \: \frac{d\gamma}{dt} 
\: \gamma^{-1} \: U(t), \\
\frac{dV}{dt} \: & = \: - \: V(t) \:  \frac{d\gamma}{dt} 
\: \gamma^{-1} \: (I - P_{\overline{\Image(d)}}(t)). \notag
\end{align}
The solution of the first equation in (\ref{4.27}), for example,
is given by
\begin{align} \label{4.28}
U(t) \: = \: & I \: + \: \int_0^t 
P_{\overline{\Image(d)}}(s) \: \frac{d\gamma}{ds} 
\: \gamma^{-1}(s) \: ds \: + \\
& \int_{t \ge s_1 \ge s_2 \ge 0}
P_{\overline{\Image(d)}}(s_1) \: \frac{d\gamma}{ds_1} 
\: \gamma^{-1}(s_1) \: P_{\overline{\Image(d)}}(s_2) \: \frac{d\gamma}{ds_2} 
\: \gamma^{-1}(s_2)
\: ds_1 \: ds_2 \: + \: \ldots \notag
\end{align}
The series in (\ref{4.28}) is convergent because 
$\frac{d\gamma}{ds} 
\: \gamma^{-1}(s)$ is uniformly bounded for $s \in [0, t]$.
One can write a similar series for $U(t)^{-1}$, showing that
$U(t)$ is invertible.

Hence $\left[ P_{\overline{\Image(d)}}(t), a \right]$ is
compact for all $t \in [0, 1]$. A similar argument 
shows that $\left[ P_{\overline{\Image(d^*)}}(t), a \right]$ is
compact for all $t \in [0, 1]$. Thus 
$\left[ F(t), a \right]$ is
compact for all $t \in [0, 1]$.

Next, $a (F(t)^2 - 1) \:  = \: -  \: a P_{\mathcal H}$, and
\begin{align} \label{4.29}
\frac{d}{dt} \: a P_{\mathcal H} \: & = \:
a \: \left( - \:
P_{\mathcal H} \: 
\frac{d\gamma}{dt} \: \gamma^{-1} \: P_{\overline{\Image(d^*)}}
\: + \: P_{\overline{\Image(d)}}  \: 
\frac{d\gamma}{dt} \: \gamma^{-1} \: P_{\mathcal H} \right) \\
& = - \: \: a \: 
P_{\mathcal H} \: 
\frac{d\gamma}{dt} \: \gamma^{-1} \: P_{\overline{\Image(d^*)}}
\: + \: [a,  \: P_{\overline{\Image(d)}}]  \: 
\frac{d\gamma}{dt} \: \gamma^{-1} \: P_{\mathcal H} 
\: + \: P_{\overline{\Image(d)}}  \: 
\frac{d\gamma}{dt} \: \gamma^{-1} \: a \: P_{\mathcal H}. \notag
\end{align}
Putting $M(0) \: = \: N(0) \: = \: I$ and solving
\begin{align} \label{4.30}
\frac{dM}{dt} \: & = \: - \: M(t) \: P_{\overline{\Image(d)}}(t) \: 
\frac{d\gamma}{dt} 
\: \gamma^{-1}, \\
\frac{dN}{dt} \: & = \: \frac{d\gamma}{dt} 
\: \gamma^{-1} \: P_{\overline{\Image(d^*)}}(t) \: N(t), \notag
\end{align}
we can write
\begin{equation} \label{4.31}
M(t) \: a \: P_{\mathcal H}(t) \: N(t) \: - \: a \: P_{\mathcal H}(0) \: = \:
\int_{0}^t M(s) \: [a,  \: P_{\overline{\Image(d)}}(s)]  \: 
\frac{d\gamma}{ds} \: \gamma^{-1} \: P_{\mathcal H}(s) \: N(s) \: ds.
\end{equation}
As $M(t)$ and $N(t)$ are invertible and 
$a P_{\mathcal H}(0)$ is compact, it follows that
$a P_{\mathcal H}(t)$ is compact for all $t \in [0, 1]$.
\end{proof}

\begin{corollary} \label{4.32}
If $\phi : X_1 \rightarrow X_2$ is an orientation-preserving
$K$-quasiconformal homeomorphism, for some $K < \infty$,
and $X_1$ satisfies Assumption \ref{4.13}, then
$(H_2, \gamma_2, F_2)$ defines a cycle for
$\KKK_{2k}(C_0(X_2); \C)$.
\end{corollary}
\begin{proof}
This follows from the proof of Proposition \ref{4.23}.
\end{proof}

\begin{corollary} \label{4.33}
\cite[Theorem 1.1]{Hilsum (1985)},
\cite[p. 678]{Connes-Sullivan-Teleman (1994)}
If $\phi \: : \: X_1 \rightarrow X_2$ is an orientation-preserving 
homeomorphism between compact oriented smooth manifolds then the
pushforward of the
signature class of $X_1$ coincides with
the signature class of $X_2$, in 
$\KKK_{2k}(C(X_2); \C)$.
\end{corollary}
\begin{proof}
If $\dim(X) \neq 4$ then there is an orientation-preserving 
quasiconformal 
homeomorphism from $X_1$ to $X_2$ that is isotopic to $\phi$
\cite{Sullivan (1979b)}, and the corollary follows from Proposition 
\ref{4.23}.
If $\dim(X) = 4$ then one can instead consider $X \times S^2$.
\end{proof}
\noindent
{\bf Remark : } If $X^\prime \: = \: X - Z$, where $Z$ has Hausdorff 
dimension at most $2k-2$, then the cycle $(H, \gamma, F)$
for $\KKK_{2k}(C_0(X); \C)$ also defines a signature cycle
for $\KKK_{2k}(C_0(X^\prime); \C)$. This is because the triple 
$(H, \gamma, F)$ is
the same as the corresponding triple for $X^\prime$, and an element
$a \in C_0(X^\prime)$ extends by zero to an element of $C_0(X)$. 
For example, writing $\R^{2k} \: = \: S^{2k} \: - \: \pt$, we obtain
a cycle $(H, \gamma, F)$ for $\KKK_{2k}(C_0(\R^{2k}); \C)$. \\

\subsection{When the limit set is the entire sphere, even-dimensional}
\label{Subsubsection 4.2.4}

In this subsection we use $F$ to construct an equivariant K-cycle for
$C(\Lambda)$ when $\Lambda \: = \: S^{2k}$.

Suppose that $\Lambda \: = \: S^{2k}$. The triple
$(H, \gamma, F)$ of Section \ref{Subsubsection 4.2.1} is $\Gamma$-equivariant
and so gives a cycle for a class $[(H, \gamma, F)] \in
\KKK_{2k}^\Gamma(C(S^{2k}); \C)$. As the nonequivariant K-homology
class represented by
$(H, \gamma, F)$ is the signature class,
it follows from the discussion of Section
\ref{4.1} that $[(H, \gamma, F)]$ is a nontorsion element of
$\KKK_{2k}^\Gamma(C(S^{2k}); \C)$

\section{From even cycles to odd cycles} \label{Subsection 4.3}

In this section we consider a manifold
$X$ as in Section \ref{Subsection 4.2} 
equipped with a partial compactification $\overline{X}$.
Putting $\partial \overline{X} = \overline{X} - X$, we give a sufficient
condition for the triple
$(H, \gamma, F)$ to extend to a cycle for
$\KKK_{2k}(C_0(\overline{X}), C(\partial \overline{X}); \C)$.
We then consider
the boundary map 
$\KKK_{2k}(C_0(\overline{X}), C(\partial \overline{X}); \C)
\rightarrow \KKK_{2k-1}(C(\partial \overline{X}); \C)$.
We describe the image of the cycle $(H, \gamma, F)$ as an element of
$\Ext(C(\partial \overline{X}))$. If $\partial
\overline{X}$ is a manifold then the relevant
Hilbert space turns out to be
the exact $k$-forms on 
$\partial \overline{X}$ of a certain regularity.
In the special case when $\partial \overline{X} = S^{2k-1}$, we show that
the Hilbert space of such $H^{-1/2}$-regular forms is
M\"obius-invariant, along with 
the Ext element.

A second technical assumption arises in this section, which will
again be satisfied in the cases that are relevant for limit sets.

\subsection{A relative K-cycle}
\label{Subsubsection 4.3.1}

In this subsection we start with 
a partial compactification $\overline{X}$ of $X$. Applying
the boundary map to
the K-cycle $(H, \gamma, F)$ for $C_0(X)$ gives a class in
$\KKK_{2k-1}(C(\partial \overline{X}); \C)$. We show the
compatibility of this map with quasiconformal homeomorphisms.
If $X$ is the domain of discontinuity $\Omega$ for $\Gamma$ then
we discuss the twisting of this construction by the pullback of a
vector bundle on $\Omega/\Gamma$.

Let $\overline{X}$ be a locally compact Hausdorff space that contains
$X$ as an open dense subset. 
Put $\partial \overline{X} \: = \: \overline{X} - X$, which we assume
to be compact. There is
a short exact sequence of $C^*$-algebras
\begin{equation} \label{4.34}
0 \longrightarrow C_0(X) \longrightarrow C_0(\overline{X}) 
\longrightarrow C(\partial \overline{X}) 
\longrightarrow 0.
\end{equation}
From
\cite[Theorem (14.24)]{Baum-Douglas (1991)}, \cite{Kasparov (1991)} or
\cite[Theorem 5.4.5]{Higson-Roe (2000)},
there is an isomorphism $\KKK_{2k}(C_0(X); \C) \cong
\KKK_{2k}(C_0(\overline{X}), C(\partial \overline{X}); \C)$.
Furthermore, there is a boundary map $\partial \: : \:
\KKK_{2k}(C_0(\overline{X}), C(\partial \overline{X}); \C) \rightarrow
\KKK_{2k-1}(C(\partial \overline{X}); \C)$.

Let $e \in M_N(C^\infty(X))$ be a projection. If $(H, \gamma, F)$ is
a K-cycle for $C_0(X)$ then there is a
new K-cycle $(H_e, \gamma_e, F_e)$, where
$H_e \: = \: H^N e$, $\gamma_e \: = \: e \gamma e$ and
$F_e \: = \: e F e$. In this way, we obtain a map
$\KK^0(X) \rightarrow \KKK_{2k}(C_0(X); \C)$. Composing with the
boundary map gives a map
$\KK^0(X) \rightarrow \KKK_{2k}(C_0(X); \C) \stackrel{\partial}{\rightarrow}
\KKK_{2k-1}(C(\partial \overline{X}); \C)$.

In this paragraph
we take $X \: = \Omega \neq \emptyset$ and
$\overline{X} \: = \: S^{2k}$, so $\partial \overline{X} \: = \: \Lambda$.
If $X$ satisfies Assumption \ref{4.13} then we have the K-cycle
$(H, \gamma, F)$ of Section \ref{Subsubsection 4.2.2}.
Let $p \in 
M_N(C^\infty(\Omega/\Gamma))$ be a projection. If $\pi \: : \:
\Omega \rightarrow \Omega/\Gamma$ is the quotient map then $e \: = \: \pi^* p$
is a projection in $M_N(C^\infty(\Omega))$. 
Applying the preceding construction
and taking into account the $\Gamma$-equivariance, we obtain maps
\begin{equation} \label{4.35}
\KK^0(\Omega) \rightarrow \KKK_{2k-1}(C(\Lambda); \C)
\end{equation}
and
\begin{equation} \label{4.36}
\KK^0(\Omega/\Gamma) \rightarrow \KKK_{2k-1}^\Gamma(C(\Lambda); \C).
\end{equation}
With reference to Proposition \ref{4.1},
the maps (\ref{4.35}) and (\ref{4.36}) are rationally the same as the 
connecting maps
\begin{equation} \label{4.37}
\KK^0(\Omega) \rightarrow
\KK^1(S^{2k}, \Omega) \cong \KKK_{2k-1}(C(\Lambda); \C)
\end{equation}
and
\begin{equation} \label{4.38}
\KK^0(\Omega/\Gamma) \cong \KK^0_\Gamma(\Omega) \rightarrow
\KK^1_\Gamma(S^{2k}, \Omega) \cong \KKK^\Gamma_{2k-1}(C(\Lambda); \C).
\end{equation}
We obtain a
rational instead of integral statement because the K-homology classes
defined by the signature and Dirac operator on $S^{2k}$, the latter
being the fundamental class, are only
rationally equivalent.

Returning to general $X$,
let $X^\prime$ be another manifold as in Section 
\ref{Subsubsection 4.2.1},
with partial compactification 
$\overline{X}^\prime$ and boundary $\partial \overline{X}^\prime$.
Let $\phi \: : \: \overline{X}^\prime \rightarrow \overline{X}$ be a
homeomorphism that restricts to a $K$-quasiconformal homeomorphism from
$X^\prime$ to $X$. 
By naturality, there is an isomorphism
$\left( \phi \big|_{\partial \overline{X}^\prime} \right)_* \: : \: \KKK_{2k-1}
(C(\partial \overline{X}^\prime); \C) \rightarrow 
\KKK_{2k-1}(C(\partial \overline{X}); \C)$.
Suppose that $X^\prime$ satisfies Assumption \ref{4.13}. By Proposition 
\ref{4.14}, Proposition \ref{4.16}
and Corollary \ref{4.32}, there are well-defined signature classes
$[(H^\prime, \gamma^\prime, F^\prime)] \in \KKK_{2k}(C_0(
{X^\prime}); \C) \cong \KKK_{2k}(C_0(
\overline{X^\prime}), C(\partial{\overline{X}}^\prime); \C)$ and
$[(H, \gamma, F)] \in \KKK_{2k}(C_0(
{X}); \C) \cong \KKK_{2k}(C_0(
\overline{X}), C(\partial \overline{X}); \C)$.
\begin{proposition} \label{4.40}
$\left( \phi \big|_{\partial \overline{X}^\prime} \right)_* 
(\partial [(H^\prime, \gamma^\prime, F^\prime)]) \: = \:
\partial [(H, \gamma, F)]$ in $\KKK_{2k-1}(C(\partial \overline{X}); \C)$.
\end{proposition}
\begin{proof}
There is a commutative diagram
\begin{equation} \label{4.41}
\begin{CD}
\KKK_{2k}(C_0(\overline{X}^\prime), C(\partial \overline{X}^\prime)); 
\C) @>\phi_*>> 
\KKK_{2k}(C_0(\overline{X}),
C(\partial \overline{X}); \C) \\
@V{\partial}VV @V{\partial}VV \\
\KKK_{2k-1}(C(\partial \overline{X}^\prime); \C) 
@>{\left( \phi \big|_{\partial \overline{X}^\prime} \right)_*}>> 
\KKK_{2k-1}(C(\partial \overline{X}); \C),
\end{CD}
\end{equation}
where the horizontal arrows are isomorphisms.  From Proposition \ref{4.23},
$\phi_* ([(H^\prime, \gamma^\prime, F^\prime)]) \: = \:
[(H, \gamma, F)]$.
The claim follows from the commutativity of the
diagram.
\end{proof}

\subsection{The induced structure on the boundary}
\label{Subsubsection 4.3.2}

In this subsection we consider a manifold $X$ as before with a
compactification $\overline{X}$. With an assumption on
$\overline{X}$, related to the Higson corona of $X$, we show that
the K-cycle $(H, \gamma, F)$ for $C_0(X)$ extends to a
K-cycle for $(C_0(\overline{X}), C(\partial \overline{X}))$.
We describe the Baum-Douglas
boundary map in this case.

Let $X$ be a manifold as in Section \ref{Subsubsection 4.3.1} satisfying 
Assumption
\ref{4.13}, with a partial compactification $\overline{X}$.
We recall that a relative K-cycle for the pair $(C_0(\overline{X}),
C(\partial \overline{X}))$ is given by a K-cycle $(H, \gamma, F)$ for
the ideal $C_0(X)$ so that the action of $C_0(X)$ on $H$ extends to
an action of $C_0(\overline{X})$, and for all $a \in C_0(\overline{X})$,
$[F, a] \in K(H)$. 

We wish to extend the K-cycle of Section \ref{Subsubsection 4.2.2} for
$C_0(X)$ to a K-cycle for $(C_0(\overline{X}),
C(\partial \overline{X}))$. There is an evident action of
$C_0(\overline{X})$ on $H$. We will need an additional
condition on $\overline{X}$.

\begin{assumption} \label{ass2}
With respect to a Riemannian metric on $X$ satisfying
Assumption \ref{4.13},
for each $a \in C_0(\overline{X})$, $a \big|_X$ is the norm limit of
a sequence $\{a_i\}_{i=1}^\infty$ of bounded elements of $C^\infty(X)$
satisfying $|d a_i| \in C_0(X)$.
\end{assumption}

If $\overline{X}$ is compact then Assumption \ref{ass2} is equivalent
to saying that $\partial \overline{X}$ is a quotient of the Higson corona, the
latter being defined using the given Riemannian metric on $X$. \\ \\
{\bf Example $1^\prime$ : } With reference to Example 1, 
Assumption \ref{ass2} is satisfied by an
asymptotically hyperbolic metric on $X$. \\ \\
{\bf Example $2^\prime$ : } With reference to Example 2, Assumption
\ref{ass2} is satisfied when $\overline{X} \: = \: S^2$ is the
one-point-compactification of $X$.

\begin{proposition} \label{4.39}
If Assumption \ref{ass2} is satisfied and $(H, \gamma, F)$ is the cycle for
$\KKK_{2k}(C_0(X); \C)$
from Section \ref{Subsubsection 4.2.2} 
then $(H, \gamma, F)$ is also a cycle for 
$\KKK_{2k}(C_0(\overline{X}), C(\partial \overline{X}); \C)$.
\end{proposition}
\begin{proof}
We must show that for all $a \in C_0(\overline{X})$, 
$[F, a]$ is compact. We may assume that $a \big|_X$ is
smooth and $|da| \in C_0(X)$. Then the proof of Proposition
\ref{4.16} applies.
\end{proof}

The boundary map 
$\partial \: : \: \KKK_{2k}
(C_0(\overline{X}), C(\partial \overline{X}); \C)
\rightarrow \KKK_{2k-1}(C(\partial \overline{X}); \C)$ can be explicitly
described as follows. Given $a \in C(\partial \overline{X})$, let
${a}^\prime$ be an extension of it to
$C_0(\overline{X})$. Then $P_{{\mathcal H}_\pm}  {a}^\prime
P_{{\mathcal H}_\pm}$ is an element of $B({\mathcal H}_\pm)$. The
corresponding element $[P_{{\mathcal H}_\pm} {a}^\prime 
P_{{\mathcal H}_\pm}]$ of the Calkin algebra 
$Q({\mathcal H}_\pm)$ is independent of the choice of extension and
defines an algebra homomorphism $\sigma_\pm \: : \:
C_0(\partial \overline{X}) \rightarrow
Q({\mathcal H}_\pm)$. Then $\partial [(H, \gamma, F)]$ is represented by
the Ext class $[\sigma_+] \: - \: [\sigma_-]$
\cite[Definition (4.6), Theorems (14.23) and (14.24)]{Baum-Douglas (1991)},
\cite[Remark 8.5.7]{Higson-Roe (2000)}.

\subsection{The case of a smooth manifold-with-boundary}
\label{Subsubsection 4.3.3}

In this subsection we consider the case when $\overline{X}$ is a smooth 
manifold-with-boundary. We construct a Hilbert space 
$H_{\partial \overline{X}}$
of exact $k$-forms
on $\partial \overline{X}$ as boundary values of $L^2$-harmonic $k$-forms
on $X$. There is a natural $\Z_2$-grading on the
Hilbert space coming from a diffeomorphism-invariant Hermitian
form. In the case when $\overline{X} \: = \: [0, \infty) \times 
\partial \overline{X}$, we show that the inner product on
$H_{\partial \overline{X}}$ is the $H^{-1/2}$ inner product.
 
Suppose that $\overline{X}^{2k}$ is a smooth oriented manifold-with-boundary
with compact boundary $\partial \overline{X}$.
Let $g_0$ be a smooth Riemannian metric on $\overline{X}$ and consider
the corresponding conformal class on $X$.
We assume that 
the reduced $L^2$-cohomology group $\HH^k_{(2)}(\overline{X}; \R) \: 
\cong 
\: \HH^k_{(2)}(\overline{X},
\partial \overline{X}; \R)$ vanishes.
(Note that $\HH^k_{(2)}(\overline{X}; \R)$ and
$\HH^k_{(2)}(\overline{X}, \partial \overline{X}; \R)$ have harmonic
representatives defined using
boundary conditions, and are generally much smaller than ${\mathcal H}$.)

Let $i \: : \: \partial
\overline{X} \rightarrow \overline{X}$ be the boundary inclusion.
We note that by conformal invariance, the $L^2$-harmonic $k$-forms on $X$ can
be computed using the metric $g_0$ which is smooth up to the boundary
$\partial \overline{X}$. It follows that $i^* \: : \: {\mathcal H} \rightarrow
H^{-1/2}(\partial \overline{X}; \Lambda^k)$ is well-defined
\cite[B.2.7-B.2.9]{Hormander (1985)}.

\begin{proposition} \label{4.42}
Given $\omega \in 
\Image \left( d \: : \: C^\infty(\partial \overline{X}; \Lambda^{k-1})
\rightarrow C^\infty(\partial \overline{X}; \Lambda^k) \right)$,
there is a unique ${\omega^\prime} \in {\mathcal H}$ so that
$i^* {\omega^\prime} \: = \: \omega$.
\end{proposition}
\begin{proof}
Write $\omega \: = \: d \eta$ for some 
$\eta \in 
C^\infty(\partial \overline{X}; \Lambda^{k-1})$.
Let ${\eta^\prime} \in C^\infty_c(\overline{X}; \Lambda^{k-1})$
satisfy $i^* {\eta^\prime} \: = \: \eta$. Let $G$ be the
Green's operator for the Laplacian on $\overline{X}$, as defined
using $g_0$, with
relative boundary conditions. In particular,
$i^* \circ G \: = \: 0$. If ${\omega^\prime}$ exists
then it satisfies $d({\omega^\prime} \: - \: d{\eta^\prime}) \: = 
\: 0$, $d^*({\omega^\prime} \: - \: d{\eta^\prime}) \: = 
\: - \: d^* d {\eta^\prime}$ and $i^* ({\omega^\prime} \: - \: 
d{\eta^\prime}) \: = \: 0$. These equations would imply
$\triangle ({\omega^\prime} \: - \: d{\eta^\prime}) \: = \:
- \: d d^* d  {\eta^\prime}$, which has the solution
${\omega^\prime} \: - \: d{\eta^\prime} \: = \: - \: G d d^* d 
{\eta^\prime}$. This motivates putting
${\omega^\prime} \: = \: d \left( {\eta^\prime} \:  - \: G d^* d 
{\eta^\prime} \right)$, which works. Note that $\omega^\prime$ is
square-integrable with respect to $g_0$, and hence lies in
$L^2(X; \Lambda^k)$.

If ${\omega_1^\prime}$ and ${\omega_2^\prime}$ both satisfy the
conclusion of the proposition then
$d({\omega_1^\prime} - {\omega_2^\prime}) \: = \:
d^*({\omega_1^\prime} - {\omega_2^\prime}) \: = \:
i^*({\omega_1^\prime} - {\omega_2^\prime}) \: = \: 0$.
The cohomology assumption then implies that ${\omega_1^\prime} = 
{\omega_2^\prime}$.
\end{proof}

\begin{definition} \label{4.43}
The Hilbert space $H_{\partial \overline{X}}$ is the completion of
$\Image \left( d \: : \: C^\infty(\partial \overline{X}; \Lambda^{k-1})
\rightarrow C^\infty(\partial \overline{X}; \Lambda^k) \right)$
with respect to the norm $\omega \rightarrow \parallel {\omega^\prime}
\parallel_{\mathcal H}$.
\end{definition}

\begin{corollary} \label{4.44}
If $i^* \: : \: \HH^k(\overline{X}; \C) \rightarrow
\HH^k(\partial \overline{X}; \C)$ is the zero map then
pullback gives an isometric isomorphism $i^* \: : \: {\mathcal H}
\rightarrow H_{\partial \overline{X}}$.
\end{corollary}
\begin{proof}
Given $\omega \in {\mathcal H}$, it represents a class $[\omega] \in
\HH^k(\overline{X})$. By assumption, $[i^* \omega]$ vanishes in 
$\HH^k(\partial \overline{X})$. Hence $i^* \omega \in \Image(d)$. The
lemma now follows from Proposition \ref{4.42}.
\end{proof}

\begin{definition} \label{4.45}
The operator $T \in B \left( H_{\partial \overline{X}} \right)$ is
given by
\begin{equation} \label{4.46}
T \omega \: = \:
\begin{cases}
\omega \: & \text{ if } \omega \in i^* {\mathcal H}_+, \\
- \: \omega \: & \text{ if } \omega \in i^* {\mathcal H}_-.
\end{cases}
\end{equation} 
\end{definition}

\begin{proposition} \label{4.47}
For all
$\omega_1, \omega_2 \in  
\Image \left( d \: : \: C^\infty(\partial \overline{X}; \Lambda^{k-1})
\rightarrow C^\infty(\partial \overline{X}; \Lambda^k) \right)$,
\begin{equation} \label{4.48}
\langle T\omega_1, \omega_2 \rangle \: = \: i^{k}
\: \int_{\partial \overline{X}}
\eta_1 \wedge \overline{\omega_2},
\end{equation}
where $\eta_1 \in C^\infty(\partial \overline{X}; \Lambda^{k-1})$ 
is an arbitrary solution of $d\eta_1 \: = \: \omega_1$.
\end{proposition}
\begin{proof}
Suppose that $\omega_1 \: = \: i^* \omega_1^\prime$ and
$\omega_2 \: = \: i^* \omega_2^\prime$, with 
$\omega_1^\prime, \omega_2^\prime \in {\mathcal H}$ being uniquely
determined.  Let ${\eta_1^\prime} \in 
C^\infty_c(\overline{X}; \Lambda^{k-1})$
satisfy $i^* {\eta_1^\prime} \: = \: \eta_1$. Then as in the proof of
Proposition \ref{4.42},
${\omega_1^\prime} \: = \: d \left( {\eta_1^\prime} \:  - \: G d^* d 
{\eta_1^\prime} \right)$. 

Suppose that $\omega_2^\prime \in {\mathcal H}_\pm$. Then
$* \omega_2^\prime \: = \: \pm \: i^{-k} \: \omega_2^\prime$ and so
\begin{align} \label{4.49}
\langle T\omega_1, \omega_2 \rangle \: & = \: 
\langle \omega_1, T\omega_2 \rangle \: = \: 
\pm \: 
\int_{\overline{X}}
\omega_1^\prime \wedge \overline{* \omega_2^\prime} \: = \:
i^{k} \: 
\int_{\overline{X}}
\omega_1^\prime \wedge \overline{\omega_2^\prime} \\ 
& = \:
i^{k} \: \int_{\overline{X}}
d \left( {\eta_1^\prime} \:  - \: G d^* d 
{\eta_1^\prime} \right) \wedge \overline{\omega_2^\prime} \: = \:
i^{k} \: \int_{\partial \overline{X}}
i^* \left( {\eta_1^\prime} \:  - \: G d^* d 
{\eta_1^\prime} \right) \wedge i^* \overline{\omega_2^\prime} \notag \\
& = \
i^{k} \: \int_{\partial \overline{X}} \eta_1 \wedge 
\overline{\omega_2}. \notag
\end{align}

To see directly that (\ref{4.48}) is independent of the choice of $\eta_1$,
suppose that 
$\eta_1$ and $\widetilde{\eta}_1$ satisfy
$d\eta_1 \: = \: d\widetilde{\eta}_1 \: = \: \omega_1$. Write
$\omega_2 \: = \: d\eta_2$. Then
\begin{equation} \label{4.50}
\int_{\partial \overline{X}} (\eta_1 \: - \: \widetilde{\eta}_1) 
\wedge \overline{\omega_2} \: = \:
\int_{\partial \overline{X}} (\eta_1 \: - \: \widetilde{\eta}_1)
\wedge d\overline{\eta_2} \: = \:
(-1)^k \: \int_{\partial \overline{X}} d(\eta_1 \: - \: \widetilde{\eta}_1) 
\wedge \overline{\eta_2} \: = \: 0.
\end{equation}
\end{proof}

\begin{proposition} \label{4.51}
Let $\partial \overline{X}$ be a closed oriented $(2k-1)$-dimensional
Riemannian manifold.
If $\overline{X} \: = \: [0, \infty) \times \partial \overline{X}$
then 
\begin{equation} \label{4.52}
H_{\partial \overline{X}} \: = \: 
\Image \left( d \: : \: H^{1/2}(\partial \overline{X}; \Lambda^{k-1})
\rightarrow H^{-1/2}(\partial \overline{X}; \Lambda^k) \right).
\end{equation}
\end{proposition}
\begin{proof}
The K\"unneth formula for reduced $L^2$-cohomology, along with the
fact that $[0, \infty)$ has vanishing absolute and relative
reduced $L^2$-cohomology,
implies that $\overline{X}$ has vanishing absolute and relative
reduced $L^2$-cohomology. Hence the hypotheses of Proposition
\ref{4.42} are satisfied.

If $p \: : \: \overline{X} \rightarrow \partial \overline{X}$ is
projection and $\omega \in C^\infty(\partial \overline{X}; \Lambda^k)$ 
then we will abuse notation to also write $\omega$ for $p^* \omega$.
Let $\widehat{d}$ be the exterior derivative on $\partial \overline{X}$
and let $\widehat{*}$ be the Hodge duality operator on 
$\partial \overline{X}$. Let $t$ be the coordinate on $[0, \infty)$. Then
\begin{equation} \label{4.53}
|\omega|^2 \: d\vol_{\overline{X}} \: = \:
\omega \wedge \widehat{*} \omega \wedge dt \: = \:
(-1)^{k-1} \: \omega \wedge dt \wedge \widehat{*} \omega.
\end{equation}
Hence $* \omega \: = \: (-1)^{k-1} \: dt \wedge
\widehat{*} \omega$.

Suppose that $\omega \in C^\infty(\partial \overline{X}; \Lambda^k)$ 
satisfies $\widehat{d} \: \omega \: = \: 0$ and
\begin{equation} \label{4.54}
(-i)^k \: \widehat{d} \: \widehat{*}
\omega \: = \: \lambda \: \omega
\end{equation}
with $\lambda \in \R$. If $\lambda \: > \: 0$ then
$e^{- \lambda t} \: \left( \omega \: - \: (-i)^k \: dt \:
\wedge \: \widehat{*} \omega \right) \in H_+$ and
\begin{equation} \label{4.55}
d \left( e^{- \lambda t} \: \left( \omega \: - \: (-i)^k \: dt \:
\wedge \: \widehat{*} \omega \right) \right) \: = \: 0.
\end{equation}
From the self-duality of
$e^{- \lambda t} \: \left( \omega \: - \: (-i)^k \: dt \:
\wedge \: \widehat{*} \omega \right)$, we also have
\begin{equation} \label{4.56}
d^* \left( e^{- \lambda t} \: \left( \omega \: - \: (-i)^k \: dt \:
\wedge \: \widehat{*} \omega \right) \right) \: = \: 0.
\end{equation}

Thus $\omega \in i^* {\mathcal H}_+$. Furthermore, from (\ref{4.54}),
\begin{equation} \label{4.57}
\widehat{d} \left( \frac{1}{\lambda} \: (-i)^k \: \widehat{*}
\omega \right) \: = \: \omega. 
\end{equation}
Then from Proposition \ref{4.47},
\begin{align} \label{4.58}
\langle \omega, \omega \rangle \: & = \: i^k \: \int_{\partial
\overline{X}} \left( \frac{1}{\lambda} \: (-i)^k \: \widehat{*}
\omega \right) \wedge \overline{\omega} \: = \: 
\frac{1}{\lambda} \:  \int_{\partial
\overline{X}} \widehat{*}
{\omega} \: \wedge \: \overline{\omega} \: = \:
\frac{1}{\lambda} \:  \int_{\partial
\overline{X}} \overline{\omega} \: \wedge \: \widehat{*}
\omega \\
& = \: \overline{\frac{1}{\lambda} \:  \int_{\partial
\overline{X}} \omega \: \wedge \: \overline{\widehat{*}
\omega} } \: = \:
\frac{1}{\lambda} \:  \int_{\partial
\overline{X}} \omega \: \wedge \: \overline{\widehat{*}
\omega}. \notag
\end{align}
If $\lambda \: < \: 0$ then
$e^{\lambda t} \: \left( \omega \: + \: (-i)^k \: dt \:
\wedge \: \widehat{*} \omega \right) \in H_-$ and
\begin{equation} \label{4.59}
d \left( e^{- \lambda t} \: \left( \omega \: + \: (-i)^k \: dt \:
\wedge \: \widehat{*} \omega \right) \right) \: = \: 0,
\end{equation}
so $\omega \in i^* {\mathcal H}_-$. A similar calculation gives
$\langle \omega, \omega \rangle \: = \: - \: 
\frac{1}{\lambda} \:  \int_{\partial
\overline{X}} \omega \: \wedge \: \overline{\widehat{*}
\omega}$. Thus in either case,
\begin{equation} \label{4.60}
\langle \omega, \omega \rangle \: = \: 
\frac{1}{|\lambda|} \:  \int_{\partial
\overline{X}} \omega \: \wedge \: \overline{\widehat{*}
\omega}.
\end{equation}
As the closure of
$\Image \left( d \: : \: C^\infty(\partial \overline{X}; \Lambda^{k-1})
\rightarrow C^\infty(\partial \overline{X}; \Lambda^k) \right)$ has an
orthonormal basis
given by such eigenforms, the proposition follows.
\end{proof}

\subsection{M\"obius-invariant analysis on odd-dimensional spheres}
\label{Subsubsection 4.3.4}

In this subsection we specialize the previous section to the case
$X \: = \: B^{2k}$. We show that the Hilbert space
$H_{\partial \overline{X}}$ is the 
$H^{-1/2}$ space of exact $k$-forms on $S^{2k-1}$. 
We show that M\"obius transformations of $S^{2k-1}$ act by isometries on
$H_{\partial \overline{X}}$, and quasiconformal homeomorphisms of
$S^{2k-1}$ act boundedly on
$H_{\partial \overline{X}}$.

Take $X \: = \: H^{2k}$, the upper hemisphere in $S^{2k}$,
and $\overline{X} \: = \: \overline{H^{2k}}$. 
Then 
$\HH^k_{(2)}(\overline{X}; \R) \: 
= \: \HH^k_{(2)}(\overline{X},
\partial \overline{X}; \R) \: = \: 0$ and
$i^* \: : \: \HH^k(\overline{X}; \C) \rightarrow
\HH^k(\partial \overline{X}; \C)$ is the zero map, so we can
apply
Proposition \ref{4.42} and Corollary \ref{4.44}.

\begin{corollary} \label{4.61} (c.f. \cite[Proposition 3.2]{Chen (1996)})
The group $\Isom^+(H^{2k})$ acts isometrically on
\begin{equation} \label{4.62}
H_{S^{2k-1}} \: = \: 
\Image \left( d \: : \: H^{1/2}(S^{2k-1}; \Lambda^{k-1})
\rightarrow H^{-1/2}(S^{2k-1}; \Lambda^k) \right)
\end{equation}
preserving $T$.
\end{corollary}
\begin{proof}
If $x_0 \in H^{2k}$ is a basepoint then
$\overline{H^{2k}} - x_0$ is conformally equivalent to 
$[0, \infty) \times S^{2k-1}$. The same calculations as in the
proof of Proposition \ref{4.51} show that 
\begin{equation} \label{4.63}
H_{S^{2k-1}} \: = \: 
\Image \left( d \: : \: H^{1/2}(S^{2k-1}; \Lambda^{k-1})
\rightarrow H^{-1/2}(S^{2k-1}; \Lambda^k) \right).
\end{equation}
As $\Isom^+(H^{2k})$ acts isometrically on ${\mathcal H}$, 
it acts isometrically on $H_{S^{2k-1}}$. The Hermitian
form (\ref{4.48}) is preserved by all orientation-preserving
diffeomorphisms of $\partial \overline{X}$.
\end{proof}

\begin{corollary} \label{4.64}
The group $\Isom^+(H^{2k})$ acts isometrically on
$H^{1/2}(S^{2k-1}; \Lambda^{k-1})/\Ker(d)$,
preserving the Hermitian form
\begin{equation} \label{4.65}
S(\omega_1, \omega_2) \: = \: i^k \int_{S^{2k-1}} \omega_1 
\wedge d\overline{\omega_2}.
\end{equation}
\end{corollary}
\begin{proof}
The dual space to
$ \Image \left( d \: : \: H^{1/2}(S^{2k-1}; \Lambda^{k-1})
\rightarrow H^{-1/2}(S^{2k-1}; \Lambda^k) \right)$ is
$H^{1/2}(S^{2k-1}; \Lambda^{k-1})/\Ker(d)$, which inherits an
isometric action of $\Isom^+(H^{2k})$.
The inner product on 
$H^{1/2}(S^{2k-1}; \Lambda^{k-1})/\Ker(d)$ is given by
$\omega \rightarrow \langle d\omega, G^{1/2} d\omega 
\rangle_{L^2}$.  The Hermitian form
$S$ is preserved because of its diffeomorphism invariance.
\end{proof}

We do not claim that the inner product on
$H^{1/2}(S^{2k-1}; \Lambda^{k-1})/\Ker(d)$ is
conformally invariant, i.e. invariant with respect to a conformal
change of the metric.

We remark that in the case $k=2$, $S(\omega, \omega)$ can be identified
(up to a sign) with the helicity, or asymptotic self-linking number,
of a vector field $\xi$ satisfying $i_{\xi} \: d\vol \: = \: d\omega$
\cite[Definition III.1.14, Theorem II.4.4]{Arnold-Khesin (1998)}.

\begin{proposition} \label{4.66}
An orientation-preserving quasiconformal homeomorphism
$\phi \: : \: S^{2k-1} \rightarrow S^{2k-1}$ acts boundedly
by pullback on $H^{1/2}(S^{2k-1}; \Lambda^{k-1})/\Ker(d)$,
preserving the Hermitian form $S$.
\end{proposition}
\begin{proof}
The method of proof is that of \cite[Corollary 3.2]{Nag-Sullivan (1995)},
which proves the proposition in the 
(quasisymmetric) case $k = 1$. By composing
$\phi$ with a M\"obius transformation, we may assume that $\phi$
has a fixed point $x_\infty \in S^{2k-1}$. Performing a linear
fractional transformation to send $x_\infty$ to infinity, we may
replace $S^{2k-1}$ by $\R^{2k-1}$. 
Given $\omega \in H^{1/2}(\R^{2k-1}; \Lambda^{k-1})/\Ker(d)$, 
consider its extensions $\omega^\prime \in 
H^1(\R^{2k}_+; \Lambda^{k-1})/\Ker(d)$.
Then 
\begin{equation} \label{4.67}
\parallel \omega \parallel \: = \: \inf_{\omega^\prime \: : \:
i^* \omega^\prime \: = \: \omega}
\parallel d\omega^\prime \parallel_{L^2}.
\end{equation}

There is an extension $\phi^\prime$ of $\phi$ to a $K$-quasiconformal
homeomorphism of $\R^{2k}_+$, for some $K < \infty$
\cite{Tukia-Vaisala (1982)}. The
proposition now follows from the fact that $\phi^\prime$ acts 
boundedly by pullback on $L^2(\R^{2k}_+; \Lambda^k)$.
\end{proof}

\subsection{The boundary signature operator as an Ext class}
\label{Subsubsection 4.3.5}

With $\overline{X}$ as in Section \ref{Subsubsection 4.3.3},
we show that the image of the cycle
$(H, \gamma, F)$ under the Baum-Douglas boundary map can be described
intrinsically in terms of $\partial \overline{X}$.
It is given by certain homomorphisms
from $C(\partial \overline{X})$ to the Calkin algebra of
$H_{\partial \overline{X}}$. If $\partial \overline{X} \: = \:
S^{2k-1}$ then we show that the homomorphisms are equivariant
with respect to M\"obius transformations of $S^{2k-1}$.

Suppose that $\overline{X}$ is a partial compactification as in
Section \ref{Subsubsection 4.3.3}, satisfying Assumption
\ref{ass2} and the hypothesis of Corollary \ref{4.44}.
With reference to Definition \ref{4.45},
there is a $\Z_2$-grading
$H_{\partial \overline{X}} \: = \:
H_{\partial \overline{X}, +} \: \oplus \:
H_{\partial \overline{X}, -}$ coming from $T$.
We put a smooth Riemannian metric $g_0$ on the manifold-with-boundary
$\overline{X}$ in the given conformal class.
We define $H^{-1/2}(\partial \overline{X}; \Lambda^k)$ using
the induced metric on $\partial \overline{X}$.
Let $P_{H_{\partial \overline{X}, \pm}}$ denote orthogonal
projection from $H^{-1/2}(\partial \overline{X}; \Lambda^k)$ to
$H_{\partial \overline{X}, \pm}$. From elliptic theory,
for all $a \in C(\partial \overline{X})$, 
$\left[ 
P_{H_{\partial \overline{X}, \pm}}, (1 + \triangle)^{1/4} \: a \:
(1 + \triangle)^{-1/4} \right]$ is compact.
Hence one obtains homomorphisms 
$\tau_\pm \: : \: C(\partial \overline{X}) \rightarrow 
Q \left( H_{\partial \overline{X}, \pm} \right)$ by
$\tau_\pm(a) \: = \: \left[ P_{H_{\partial \overline{X}, \pm}}
\: (1 + \triangle)^{1/4} \: a \:
(1 + \triangle)^{-1/4} \:
P_{H_{\partial \overline{X}, \pm}} \right]$.

\begin{proposition} \label{4.68}
$\partial [(H, \gamma, F)]$ equals $[\tau_+] - [\tau_-]$
in $\Ext(C(\partial \overline{X})) \cong \KKK_{2k-1}(C(\partial \overline{X}); \C)$.
\end{proposition}
\begin{proof}
We wish to show that $[\sigma_\pm] \: = \: [\tau_\pm]$. 
The method of proof is similar to that of
\cite[Proposition 4.3]{Baum-Douglas-Taylor (1989)}.
The subspace $H_{\partial \overline{X}}$ of 
$H^{-1/2}(\partial \overline{X}; \Lambda^k)$ has an induced
inner product that is boundedly equivalent to the inner product
of Definition \ref{4.43}.
To prove the proposition, it is sufficient
to use the new inner product on $H_{\partial \overline{X}}$.
Suppose first
that $a \in C^\infty(\partial \overline{X})$.
We will show that $[\sigma_\pm](a)$ equals the class of
$\left[ P_{H_{\partial \overline{X}, \pm}}
\: a \:
P_{H_{\partial \overline{X}, \pm}} \right]$ in 
$Q \left( H_{\partial \overline{X}, \pm} \right)$. From elliptic
theory, this in turn equals the class of
$\left[ P_{H_{\partial \overline{X}, \pm}}
\: (1 + \triangle)^{1/4} \: a \:
(1 + \triangle)^{-1/4} \:
P_{H_{\partial \overline{X}, \pm}} \right]$.

Let $a^\prime \in 
C^\infty_c(\overline{X})$ be an extension of $a$.
Using the isomorphism $i^* \: : \: {\mathcal H} \rightarrow
H_{\partial \overline{X}}$, it suffices to
show that $i^* P_{\mathcal H} a^\prime \: - \:
P_{H_{\partial \overline{X}}} a i^*$ is compact from
${\mathcal H}$ to $H_{\partial \overline{X}}$.
As $i^* a^\prime P_{\mathcal H} \: - \:
a P_{H_{\partial \overline{X}}} i^*$ vanishes on ${\mathcal H}$,
it suffices to show that
$i^* [P_{\mathcal H}, a^\prime] \: - \:
[P_{H_{\partial \overline{X}}}, a] i^*$ is compact.

As $P_{H_{\partial \overline{X}}}$ is a zeroth order
pseudodifferential operator, $[P_{H_{\partial \overline{X}}}, a]$ is
compact on $H^{-1/2}(\partial \overline{X}; \Lambda^k)$, so
$[P_{H_{\partial \overline{X}}}, a] i^*$ is compact from
${\mathcal H}$ to $H^{-1/2}(\partial \overline{X}; \Lambda^k)$.

From Proposition
\ref{4.39}, $[P_{\mathcal H}, a^\prime]$ is compact from
$L^2({X}; \Lambda^k)$ to $L^2({X}; \Lambda^k)$.
Let $D$ be the operator $d \: + \: d^*$ on $X$, where $d^*$ is defined
using $g_0$. Its maximal domain is $\Dom(D_{max}) \: = \:
\{ \omega \in L^2(X; \Lambda^*) \: : \: (d \: + \: d^*)\: \omega \in
L^2(X; \Lambda^*)\}$. Clearly ${\mathcal H} \subset \Dom(D_{max})$.
Applying \cite[Lemma 3.2]{Baum-Douglas-Taylor (1989)}, we conclude
that $i^* [P_{\mathcal H}, a^\prime]$ is compact from
${\mathcal H}$ to $H^{-1/2}(\partial \overline{X}; \Lambda^k)$.

If $a$ is merely continuous then multiplication by $a$ may not
be defined on $H^{-1/2}(\partial \overline{X}; \Lambda^k)$. However,
the operator $(1 + \triangle)^{1/4} \: a \:
(1 + \triangle)^{-1/4}$ is well-defined and gives a homomorphism
$C(\partial \overline{X}) \rightarrow B(H^{-1/2}(\partial \overline{X}; 
\Lambda^k))$. The proposition now follows from the norm density of
$C^\infty(\partial \overline{X})$ in 
$C(\partial \overline{X})$.
\end{proof}

Taking ${X} \subset S^{2k}$ to be the upper hemisphere
$H^{2k}$,
it follows that $[\tau_+] - [\tau_-] \in 
\Ext(C(S^{2k-1})) \cong \KKK_{2k-1}(C(S^{2k-1}); \C)$ is
the signature class of $S^{2k-1}$.

\begin{corollary} \label{4.69}
The map $\tau_\pm \: : \: C(S^{2k-1}) \rightarrow Q(H_{S^{2k-1}, \pm})$ is 
$\Isom^+(H^{2k})$-equivariant.
\end{corollary}
\begin{proof}
This follows from the fact that the proof of Proposition \ref{4.68} is
essentially
$\Isom^+(H^{2k})$-equivariant.  We give an alternative direct argument.

The group $\Isom^+(H^{2k})$ acts on 
$H^{-1/2}(S^{2k-1}; \Lambda^k)$ through its action on 
$S^{2k-1}$, although not isometrically.  For $g \in \Isom^+(H^{2k})$,
we have
$g \: P_{H_{S^{2k-1}}} \: = \:
P_{H_{S^{2k-1}}} \:  g \: P_{H_{S^{2k-1}}}$.
Then $\Isom^+(H^{2k})$ acts by automorphisms on 
$B \left( H_{S^{2k-1}} \right)$, with 
$g \in \Isom^+(H^{2k})$ sending $T \in 
B \left( H_{S^{2k-1}} \right)$ to
$P_{H_{S^{2k-1}}} g T
g^{-1} P_{H_{S^{2k-1}}} \: = \:
P_{H_{S^{2k-1}}} g P_{H_{S^{2k-1}}} T P_{H_{S^{2k-1}}} 
g^{-1} P_{H_{S^{2k-1}}}$. There is an induced action on $Q
\left( H_{S^{2k-1}} \right)$.

Suppose that $a \in C^\infty(S^{2k-1})$ and $g \in \Isom^+(H^{2k})$.
Then 
\begin{align} \label{4.70}
P_{H_{S^{2k-1}}} g a g^{-1} P_{H_{S^{2k-1}}} \: & = \:
P_{H_{S^{2k-1}}} g a P_{H_{S^{2k-1}}}
g^{-1} P_{H_{S^{2k-1}}} \: = \:
P_{H_{S^{2k-1}}} g a P_{H_{S^{2k-1}}}^2
g^{-1} P_{H_{S^{2k-1}}} \\
& = \:
P_{H_{S^{2k-1}}} g \left( a P_{H_{S^{2k-1}}} -
P_{H_{S^{2k-1}}} a \right) P_{H_{S^{2k-1}}}
g^{-1} P_{H_{S^{2k-1}}} \: + \notag \\
& \: \: \: \: \: \:
P_{H_{S^{2k-1}}} g P_{H_{S^{2k-1}}} a P_{H_{S^{2k-1}}}
g^{-1} P_{H_{S^{2k-1}}}. \notag
\end{align}
From elliptic theory, $a P_{H_{S^{2k-1}}} -
P_{H_{S^{2k-1}}} a$ is compact.  It follows that the homomorphism
$C^\infty(S^{2k-1}) \rightarrow Q(H_{S^{2k-1}, \pm})$ is
$\Isom^+(H^{2k})$-equivariant. The corollary now follows by continuity.
\end{proof}

\section{Odd cycles on limit sets} \label{Subsection 4.4}

In this section we construct $\Gamma$-equivariant Ext cycles on limit
sets.  If the limit set is the entire sphere-at-infinity $S^{2k-1}$
then we use the Ext cycle of Section \ref{Subsubsection 4.3.5}. 
If the limit set is
a proper
subset of the sphere-at-infinity $S^{2k}$ then we take $X$ to be
a $\Gamma$-invariant union of connected components of the
domain-of-discontinuity $\Omega$. We apply the boundary construction
of Section \ref{Subsubsection 4.3.2} 
to get an Ext cycle on $\Lambda$. We show that
the resulting K-homology class is invariant under quasiconformal deformation.
We use Section \ref{Subsubsection 4.3.5} 
to describe an explicit Ext cycle for the
K-homology class in the quasiFuchsian case, and in the case of
an acylindrical convex-cocompact hyperbolic 3-manifold with 
incompressible boundary.

\subsection{When the limit set is the entire sphere,
odd-dimensional}
\label{Subsubsection 4.4.1}

In this subsection we 
suppose that $n \: = \: 2k-1$ and $\Lambda \: = \: S^{2k-1}$.

From Corollary \ref{4.69}, we have $\Gamma$-equivariant homomorphisms
$\tau_\pm \: : \: C(S^{2k-1}) \rightarrow Q(H_{S^{2k-1}, \pm})$.
In the nonequivariant case the difference of such homomorphisms 
defines an Ext class and hence an 
odd KK-class, as the relevant algebra $C(S^{2k-1})$ is nuclear
\cite[Corollary 5.2.11 and Theorem 8.4.3]{Higson-Roe (2000)}.
In the equivariant case an odd KK-class gives rise to a
$\Gamma$-equivariant Ext class, but the converse is not automatic
(see \cite{Thomsen (2000)}).
However, it is true in our case, where the relevant KK-class is the image
of the signature class of $B^{2k}$ under the maps
$\KKK^\Gamma_{2k}(C_0(B^{2k}); \C) \cong
\KKK^\Gamma_{2k}(C(\overline{B^{2k}}), C(S^{2k-1}); \C) 
\stackrel{\partial}{\rightarrow} 
\KKK^\Gamma_{2k-1}(C(S^{2k-1}); \C)$. From the discussion of
Section \ref{Subsection 4.1},
this is a nontorsion class.

\subsection{Quasiconformal invariance II}
\label{Subsubsection 4.4.2}

In this subsection we 
take $X$ to be
a $\Gamma$-invariant union of connected components of the
domain-of-discontinuity $\Omega$. We give sufficient conditions for
Assumption \ref{ass2} to be satisfied. We show that the K-homology
class arising from the boundary construction
of Section \ref{Subsubsection 4.3.2} 
is invariant under quasiconformal deformation.

Let $\Gamma^\prime$ be a discrete torsion-free subgroup of $\Isom^+(H^{2k+1})$,
with limit set $\Lambda^\prime$ and domain of discontinuity $X^\prime \: = \:
\Omega^\prime$.
We take the compactification $\overline{X}^\prime \: = \: S^{2k}$.

\begin{proposition} \label{4.71}
1. If $\Lambda^\prime =
S^{2k-l}$ and $l \neq 2$ then the compactification satisfies
Assumption \ref{ass2}. \\
2. If $\Gamma^\prime$ is convex-cocompact but not cocompact,
and the convex core has
totally geodesic boundary, then the compactification satisfies
Assumption \ref{ass2}.
\end{proposition}
\begin{proof}
1. If $\Lambda^\prime \: = \: S^{2k-l}$ then 
$\Omega^\prime$ is conformally equivalent
to $H^{2k-l+1} \times S^{l-1}$. Consider the metric on
$H^{2k-l+1} \times S^{l-1}$ that is a product of constant-curvature
metrics.  If $l$ is odd then the differential form Laplacian on 
$H^{2k-l+1}$ has a gap away from zero in its spectrum.  It follows that
Assumption \ref{4.13} is satisfied in this case.  If $l$ is even then the $p$-form
Laplacian on $H^{2k-l+1}$ is strictly positive if $p \neq
k-\frac{l}{2}, k-\frac{l}{2} + 1$. From this, the $p$-form Laplacian on
$H^{2k-l+1} \times S^{l-1}$ is strictly positive if
$p \neq
k-\frac{l}{2}, k-\frac{l}{2} + 1, k+\frac{l}{2}-1, k+\frac{l}{2}$.
It follows that the $k$-form Laplacian on 
$H^{2k-l+1} \times S^{l-1}$ is strictly positive if $l \neq 2$.
As the inclusion $\Omega^\prime \rightarrow S^{2k}$ factors through
continuous maps $\Omega^\prime \rightarrow \overline{H^{2k-l+1}} \times S^{l-1}
\rightarrow S^{2k}$, it follows that Assumption \ref{ass2} is
satisfied. \\
2. In this case $\Omega^\prime$ is a union of round balls in $S^{2k}$ 
with disjoint closures. Putting the hyperbolic metric on each of these
balls, Assumption \ref{ass2} is satisfied.
\end{proof}

There is an evident extension of Proposition \ref{4.71}.2 to the case when
rank-$2k$ cusps are allowed.

Let $\Gamma$ and $\Gamma^\prime$
be discrete torsion-free subgroups of 
$\Isom^+(H^{2k+1})$. They are said to be quasiconformally related if
there are an isomorphism $i \: : \: \Gamma^\prime
\rightarrow \Gamma$ 
and a quasiconformal homeomorphism 
$\phi \: : \: S^{2k} \rightarrow S^{2k}$ satisfying
\begin{equation} \label{4.72}
\phi \: \circ \: \gamma^\prime \: \circ \: \phi^{-1} \: = \: 
i(\gamma^\prime)
\end{equation}
for all $\gamma^\prime \in \Gamma^\prime$. It follows that the limit
sets $\Lambda^\prime$ and $\Lambda$ are related by
$\phi(\Lambda^\prime) \: = \: \Lambda$.

Let $X^\prime$ be a $\Gamma^\prime$-invariant union of connected components of
$\Omega^\prime$. Suppose that $X^\prime$ satisfies Assumption \ref{ass2}.
Then the construction described in Section \ref{Subsubsection 4.3.2}
gives $\Gamma^\prime$-equivariant homomorphisms $\sigma_\pm \: : \:
C(\Lambda^\prime) \rightarrow Q(H_{\partial \overline{X}^\prime, \pm})$.
As in the previous section, the equivariant Ext class $[\sigma_+] \: - \:
[\sigma_-]$ arises from a class in 
$\KKK_{2k-1}^\Gamma(C(\Lambda^\prime); \C)$. 

Suppose that $\Gamma$ and $\Gamma^\prime$ are quasiconformally
related.
By naturality, there is an isomorphism
$\left( \phi \big|_{\Lambda^\prime} 
\right)_* \: : \: \KKK^{\Gamma^\prime}_{2k-1}
(C(\Lambda^\prime); \C) \rightarrow 
\KKK^\Gamma_{2k-1}(C(\Lambda); \C)$.
Put $X \: = \: \phi(X^\prime)$.
Then
$\partial \overline{X^\prime} \: = \: \Lambda^\prime$ and
$\partial \overline{X} \: = \: \Lambda$.
Suppose that $X^\prime$ satisfies Assumption \ref{4.13}. By Proposition 
\ref{4.14}, Proposition \ref{4.16}
and Corollary \ref{4.32}, there are well-defined signature classes
$[(H^\prime, \gamma^\prime, F^\prime)] \in \KKK^{\Gamma^\prime}_{2k}(C(
{X^\prime}); \C) \cong \KKK^{\Gamma^\prime}_{2k}(C(
\overline{X^\prime}), C(\Lambda^\prime); \C)$ and
$[(H, \gamma, F)] \in \KKK^{\Gamma}_{2k}(C(
{X}); \C) \cong \KKK^{\Gamma}_{2k}(C(
\overline{X}), C(\Lambda); \C)$.
\begin{proposition} \label{4.73}
$\left( \phi \big|_{\Lambda^\prime} 
\right)_* (\partial [(H^\prime, \gamma^\prime, F^\prime)]) \: = \:
\partial [(H, \gamma, F)]$ in $\KKK^\Gamma_{2k-1}(C(\Lambda); \C)$.
\end{proposition}
\begin{proof}
The proof is the same as that of Proposition \ref{4.40}, extended to
the equivariant setting.
\end{proof}

Given a discrete group $G$, it follows that quasiconformally equivalent
embeddings $G \rightarrow \Isom^+(H^{n+1})$ give rise to the
same KK-class.  We note that if $\Gamma$ is a convex-cocompact
representation of $G$ then $G$ is Gromov-hyperbolic and 
$\Lambda$ is homeomorphic to $\partial G$. 
In principle the K-cycle that we have constructed
for $\KKK^\Gamma_{2k-1}(C(\Lambda); \C)$
can be expressed entirely in terms of $G$.

\subsection{Odd-dimensional quasiFuchsian manifolds}
\label{Subsubsection 4.4.3}

In this subsection we give an explicit $\Gamma$-equivariant Ext cycle
for the K-homology class in the quasiFuchsian case, as a pushforward
of the Fuchsian cycle.

Let $\Gamma^\prime$ be a discrete torsion-free subgroup of 
$\Isom^+(H^{2k})$ whose limit set is $S^{2k-1}$.
There is a natural Fuchsian embedding $\Gamma^\prime \subset \Isom^+(H^{2k+1})$.
Take $X^\prime \: = \: B^{2k}$, the upper hemisphere.
By Proposition \ref{4.71}.1, Assumption \ref{ass2} is satisfied.
A group $\Gamma \subset \Isom^+(H^{2k+1})$ 
that is quasiconformally related to $\Gamma^\prime$
is said to be a quasiFuchsian deformation of $\Gamma^\prime$.

\begin{corollary} \label{4.74}
$\partial [(H, \gamma, F)]$ is the pushforward under 
$\phi \big|_{S^{2k-1}}$ of
the signature class of
$S^{2k-1}$ in $\KKK^\Gamma_{2k-1}(C(S^{2k-1}); \C)$.
\end{corollary}
\begin{proof}
This follows from Proposition \ref{4.73}.
\end{proof}

The Ext cycle for the signature class of
$S^{2k-1}$ in $\KKK^\Gamma_{2k-1}(C(S^{2k-1}); \C)$ was
described in Section \ref{Subsubsection 4.4.1}.
Given the quasiFuchsian group $\Gamma$,
suppose that $\phi_1$ and $\phi_2$ are two quasiconformal maps
satisfying (\ref{4.72}). Then $\phi_1^{-1} \: \circ \phi_2 \big|_{S^{2k-1}}
\: : \: S^{2k-1} \rightarrow S^{2k-1}$ commutes with each element of
$\Gamma^\prime$. As the fixed points of the hyperbolic elements of
$\Gamma^\prime$ are dense in its limit set $S^{2k-1}$, it follows that
$\phi_1^{-1} \: \circ \phi_2 \big|_{S^{2k-1}} \: = \: \Id_{S^{2k-1}}$,
so $\phi_1 \big|_{S^{2k-1}} \: = \: \phi_2 \big|_{S^{2k-1}}$. 
Next, suppose that $\Gamma^{\prime \prime}$ is another Fuchsian group
such that
$H^{2k}/\Gamma^\prime$ is orientation-preserving
isometric to $H^{2k}/\Gamma^{\prime \prime}$.
Then there is some $g \in \Isom^+(H^{2k})$ so that
$g \Gamma^\prime g^{-1} \: = \: \Gamma^{\prime \prime}$. As
$g$ acts conformally on $S^{2k-1}$, we can {\em define} a conformal structure
on $\Lambda$ to be the standard conformal structure on the homeomorphic set
$\phi^{-1} (\Lambda) \: = \: S^{2k-1}$.  
This is independent of the choices made.

The upshot is that there is a $\Gamma$-equivariant
Ext cycle for the K-homology class in
$\KKK^\Gamma_{2k-1}(C(\Lambda); \C)$, given by the
pushforward of the signature Ext cycle for $S^{2k-1}$ under the
homeomorphism $\phi \big|_{S^{2k-1}} \: : \:
S^{2k-1} \rightarrow \Lambda$. From Section \ref{Subsubsection 4.4.1},
the signature Ext class for $S^{2k-1}$ is nontorsion in 
$\KKK_{2k-1}^{\Gamma^\prime}(C(S^{2k-1}); \C)$. As $\left( 
\phi \big|_{S^{2k-1}} \right)_*$ is an isomorphism, it follows that
the class in $\KKK_{2k-1}^\Gamma(C(\Lambda); \C)$ is also nontorsion.

\subsection{The case of a quasicircle}
\label{Subsubsection 4.4.4}

Applying the construction of
Section \ref{Subsubsection 4.4.3} in the case $k = 1$, we show that we
recover the
K-homology class on a quasicircle considered by Connes and Sullivan.

Suppose that $k = 1$ and $\Gamma \subset \Isom^+(H^3)$ is a 
quasiFuchsian group. Let $B^2$ be the open upper hemisphere in $S^2$
and put $X \: = \: \phi(B^2)$. If $D^2$ is the closed disk in $\C$, let 
$Z \: : \: \Int(D^2) \rightarrow X$ be a uniformization,
i.e. a holomorphic isomorphism.
The pullback $Z^* \: : \:
L^2(X; \Lambda^1) \rightarrow L^2(\Int(D^2); \Lambda^1)$ is an
isometry. Because $Z$ is a conformal diffeomorphism, $Z^*$ 
sends ${\mathcal H}_{X}$ isometrically to
${\mathcal H}_{\Int(D^2)}$. More explicitly, the elements of
${\mathcal H}_{\Int(D^2)}$ are square-integrable forms
$f_1(z) dz \: + \: \overline{f_2(z)} d\overline{z}$ on
$\Int(D^2)$, where
$f_1$ and $f_2$ are holomorphic functions on $\Int(D^2)$.

By Carath\'eodory's theorem, $Z$
extends to a homeomorphism $Z \: : \: D^2 \rightarrow
\overline{X}$ \cite[Theorem 14.19]{Rudin (1974)}.
Then
$Z^* H_{\partial \overline{X}}$ is isometric to
$\Image \left( d \: : \: H^{1/2}(S^1; \Lambda^0) \rightarrow
H^{-1/2}(S^1; \Lambda^1) \right)$, with the operator $T$ acting by
\begin{equation} \label{4.75}
T \left( e^{ik\theta} d\theta \right) \: = \:
\begin{cases}
e^{ik\theta} d\theta & \text{ if } k > 0, \\
- \: e^{ik\theta} d\theta & \text{ if } k < 0.
\end{cases}
\end{equation}
Unequivariantly, 
the homomorphisms $\sigma_{\pm} \: : \:
C(S^1) \rightarrow Q(H_{S^1, \pm})$ are essentially
the same as the standard
Toeplitz homomorphisms.

We remark that
the dual space to
$Z^* H_{\partial \overline{X}}$ is $H^{1/2}(S^1; \Lambda^0)/\C$.
The Hermitian form
$S(f_1, f_2) \: = \: \int_{S^1} f_1 \wedge d \overline{f_2}$ on
$H^{1/2}(S^1; \Lambda^0)/\C$ is the
Hermitian form of the Hilbert transform. 

Let us compare the equivariant Ext class $[\sigma_+] - [\sigma_-]$ with
that considered by Connes and Sullivan
\cite[Section IV.3.$\gamma$]{Connes (1994)}. The latter is based on the
Hilbert space $H_0 \: = \: L^2(S^1)$.  The obvious $\Gamma$-action on
$H_0$ is not unitary, but one can make it unitary by adding compensating
weights. Then there is a $\Gamma$-invariant operator $T_0$ on $H_0$,
which is essentially the Hilbert transform, and satisfies $T_0^2 \: = \: 1$.
Decomposing $H_0$ with respect to $T_0$ as $H_0 \: = \: H_{0,+} \:
\oplus H_{0,-}$, one obtains $\Gamma$-invariant homomorphisms
$\sigma_{0, \pm} \: : \: C(S^1) \rightarrow Q(H_{0,\pm})$ given by
$\sigma_{0, \pm}(f) \: = \: \frac{1 \pm T_0}{2} \: f \: \frac{1 \pm T_0}{2}$,
modulo $K(H_{0,\pm})$.

Although there is a formal similarity between $H_{S^1, \pm}$ and
$H_{0,\pm}$, they carry distinct representations of $\Gamma$. Nevertheless,
the ensuing classes in $\KK_1^\Gamma(C(S^1); \C)$ are the same.
To see this, consider the $E_2$-term $E_2^{0,0} \: = \:
H^0(\Gamma; K_1(S^1))$  in the proof of Proposition \ref{4.1}. This
term is unaffected by the differentials of the spectral sequence and
passes to the limit to give a contribution to $\KK_1^\Gamma(C(S^1); \C)$.
It corresponds to $\Gamma$-invariant elements of $K_1(S^1)$. Unequivariantly,
$[\sigma_{+}] \: - \: [\sigma_{-}] \: = \: 
[\sigma_{0,+}] \: - \: [\sigma_{0,-}]$ in $K_1(S^1)$. As both sides
are $\Gamma$-invariant, it follows that
they give rise to the same class in $\KK_1^\Gamma(C(S^1); \C)$.

We note that the main use of the Connes-Sullivan cycle is to define
certain operators on $H_0$ for which one wants to compute the trace.
As the trace is formally independent of the choice of inner product,
one can consider the same operators on $H_{S^1}$.  See the remark after
Proposition \ref{4.89} for further discussion. 

\subsection{Odd-dimensional convex-cocompact manifolds}
\label{Subsubsection 4.4.5}

In this subsection we give an explicit $\Gamma$-equivariant Ext cycle
in the case of an odd-dimensional
convex-cocompact hyperbolic
manifold whose convex core has totally geodesic boundary.  We use this to 
give an explicit cycle in the case of an arbitrary acylindrical 
convex-cocompact 
hyperbolic $3$-manifold with incompressible boundary.

Let $M^{2k+1}$ be a noncompact convex-cocompact hyperbolic manifold
with a convex core $Z \subset M$ whose boundary is totally geodesic.
Let $C$ be a boundary component of $\partial \overline{M}$. Then the 
preimage $X$ of $C$ in $\Omega$ is 
a union $\bigcup_{i=1}^\infty B_i$ of round balls in
$S^{2k}$ with disjoint closures. Put $Y_i \: = \: \partial
\overline{B_i}$. Then $\Lambda$ is the closure of
$\bigcup_{i=1}^\infty Y_i$. By Proposition \ref{4.71}.2,
Assumption \ref{ass2} is satisfied.
We now describe the Ext cycle on $\Lambda$ coming from
Section \ref{Subsubsection 4.3.2}.
From Section \ref{Subsubsection 4.3.4}, the Hilbert space
will be $H \: = \: \bigoplus_{i=1}^\infty 
\Image \left( d \: : \: H^{1/2}(Y_i; \Lambda^{k-1}) \rightarrow
H^{-1/2}(Y_i; \Lambda^{k}) \right)$. It is $\Z_2$-graded by the
operator $T$ of Definition \ref{4.45}, applied separately to each
$Y_i$. The Ext class will be $[\sigma_+] - [\sigma_-]$, where the
homomorphisms $\sigma_\pm \: : \: C(\Lambda) \rightarrow
Q(H_\pm)$ come from restricting $f \in C(\Lambda)$ to each
$Y_i$ and applying the map $\tau_\pm$ of Corollary
\ref{4.69}.

Now let $M$ be a noncompact acylindrical convex-cocompact 
hyperbolic $3$-manifold with incompressible boundary.
Let $Z$ be a compact core for $M$. 
There is a hyperbolic
$3$-manifold $M^\prime$, homeomorphic to $M$, whose convex core
has totally geodesic boundary (one applies
Thurston's hyperbolization
theorem for Haken manifolds to get an involution-invariant
hyperbolic metric on the double $DZ$).
Furthermore, it follows from
\cite[Theorem 8.1]{Marden (1974)} that the groups 
$\Gamma^\prime \: = \: \pi_1(M^\prime)$ and
$\Gamma \: = \: \pi_1(M)$ are
quasiconformally related. 
The K-homology
class on $\Lambda^\prime$ is represented by the Ext cycle of the
preceding paragraph.
From Proposition \ref{4.73}, the K-homology 
class on $\Lambda$ is represented by the pushforward of this Ext cycle
by $\phi \big|_{\Lambda^\prime}$.
From the discussion of Section \ref{Subsection 4.1}, if $\partial
\overline{M}$ has more than one connected component then one gets
nontorsion K-homology classes from this construction. Topologically,
$\Lambda$ is a Sierpinski curve.

There is an evident extension to the case when $M$ is allowed to have
rank-two cusps.

\section{From odd cycles to even cycles} \label{Subsection 4.5}

In Section \ref{Subsection 4.4} 
we considered the case when $\Lambda$ is a proper subset
of $S^{2k}$ and showed how to pass from an even K-cycle on $\Omega$ to
an Ext cycle on $\Lambda$. In this section we consider the case when 
$\Lambda$ is a proper subset
of $S^{2k-1}$. We then want to start with an odd cycle on $\Omega$ and
construct an even K-cycle on $\Lambda$.

In the closed case, the relevant Hilbert space for an Ext cycle
is the dual space
to that of Section \ref{Subsubsection 4.3.3}, 
namely $H^{1/2}(X, \Lambda^{k-1})/\Ker(d)$.
If $X$ instead has a compactification $\overline{X}$ then there are
different choices for
$H^{1/2}(X, \Lambda^{k-1})/\Ker(d)$, depending on the particular metric
(complete or incomplete) taken in the given conformal class.  This point
deserves further study.  A related problem is to develop a good notion of
a relative version of Ext and the corresponding boundary map, as mentioned 
in  \cite[p. 3]{Baum-Douglas (1991)}. Of course there is a boundary map
in odd relative K-homology \cite[Proposition 8.5.6(b)]{Higson-Roe (2000)}, but 
in our case the natural cycles are Ext cycles.
In this section we 
will just illustrate using smooth forms how to go from the odd cycle
on $X$ to an even K-cycle on $\partial X$.
We describe the resulting
K-cycle in the quasiFuchsian case, and in the case of a quasiconformal
deformation of a convex-cocompact hyperbolic manifold whose convex
core has totally
geodesic boundary.
In the case $k = 1$ we recover the K-cycle on a
Cantor set considered by Connes and Sullivan.

\subsection{The boundary map in the odd case}
\label{Subsubsection 4.5.1}

In this subsection we describe a formalism to go from the Ext cycle of
Section \ref{Subsubsection 4.3.3}, considered 
on an odd-dimensional manifold-with-boundary, to
an even K-cycle on the boundary.

Let $X^{2k-1}$ be an odd-dimensional compact oriented manifold-with-boundary.
Let $i \: : \: \partial X \rightarrow X$ be the boundary inclusion.
We write
\begin{equation} \label{4.76}
\Ker(d) \: = \: \Ker \left( d \: : \: 
C^\infty(X; \Lambda^{k-1}) \rightarrow
C^\infty(X; \Lambda^{k}) \right)
\end{equation} 
and
\begin{equation} \label{4.77}
\Ker(d)_0 \: = \: \{ \omega \in 
\Ker(d) \: : \: i^* \omega \: = \: 0\}.
\end{equation} 
The form
\begin{equation} \label{4.78}
S(\omega_1, \omega_2) \: = \: i^k \int_{X} \omega_1 
\wedge d\overline{\omega_2}
\end{equation} 
is well-defined on $C^\infty(X; \Lambda^{k-1})/\Ker(d)_0$
and satisfies
\begin{equation} \label{4.79}
S(\omega_1, \omega_2) \: - \: \overline{S(\omega_2, \omega_1)} \: = \:
- \: (-i)^k \int_{\partial X} i^* \omega_1 \wedge i^* \overline{\omega_2}.
\end{equation}

The map $i^* \: : \: C^\infty(X; \Lambda^{k-1}) \rightarrow
C^\infty(\partial X; \Lambda^{k-1})$ restricts to a map on
$\Ker(d)/\Ker(d)_0$, with image
$i^* \Ker(d) \subset C^\infty(\partial X; \Lambda^{k-1})$.

We now assume that $\partial X$ has a conformal structure. Then
we have the Hilbert space $H_{\partial X} \: = \: 
L^2(\partial X; \Lambda^{k-1})$, with $\Z_2$-grading operator
$\gamma$ as in (\ref{4.5}). From (\ref{4.79}),
\begin{equation} 
S(\omega_1, \omega_2) \: - \: \overline{S(\omega_2, \omega_1)} \: = \:
(-1)^{k+1} \: i \: \langle i^* \omega_1, \gamma \: i^* \omega_2 
\rangle_{\partial X}.
\end{equation}
This is a compatibility between the form $S$ on $X$ and the
inner product on $\partial X$.

\begin{proposition} \label{4.80}
There is an orthogonal decomposition
\begin{equation} \label{4.81}
H_{\partial X} \: = \: \overline{i^* \Ker(d)} \: \oplus \:
\gamma \: \overline{i^* \Ker(d)}.
\end{equation}
\end{proposition}
\begin{proof}
Suppose that $\omega^\prime_1, \omega^\prime_2 \in \Ker(d) \subset
C^\infty(X; \Lambda^{k-1})$. Then
\begin{equation} \label{4.82}
\int_{\partial X} \omega^\prime_1 \wedge \overline{\omega^\prime_2} \: = \:
\int_X d(\omega^\prime_1 \wedge \overline{\omega^\prime_2}) \: = \: 0.
\end{equation}
This implies that $\overline{i^* \Ker(d)}$ and 
$\gamma \: \overline{i^* \Ker(d)}$ are perpendicular.

If $\omega \: = \: d \eta$ with
$\eta \in C^\infty(\partial X; \Lambda^{k-2})$, and $\eta^\prime
\in C^\infty(X; \Lambda^{k-2})$ satisfies $i^* \eta^\prime \: = \:
\eta$, then $\omega \: = \: i^* d\eta^\prime$. Thus
$\Image(d \: : \: C^\infty(\partial X; \Lambda^{k-2}) \rightarrow
C^\infty(\partial X; \Lambda^{k-1}) )$ 
is contained in $i^* \Ker(d)$, and similarly
$\Image(d^* \: : \: C^\infty(\partial X; \Lambda^{k}) \rightarrow
C^\infty(X; \Lambda^{k-1}) )$ is contained in 
$\gamma \: i^* \Ker(d)$.

Suppose that $\omega \in H_{\partial X}$ is orthogonal to
$\overline{i^* \Ker(d)}$ and 
$\gamma \: \overline{i^* \Ker(d)}$. 
It follows that
$d \omega \: = \: d^* \omega \: = \: 0$. Without loss of generality, 
we can take $\omega$ to be real.
Let $[\omega] \in
\HH^{k-1}(\partial X; \R)$ denote the corresponding cohomology class.
From the cohomology exact sequence
\begin{equation} \label{4.83}
\ldots \rightarrow \HH^{k-1}(X; \R) \stackrel{i^*}{\rightarrow} 
\HH^{k-1}(\partial X; \R) \stackrel{(i^*)^*}{\rightarrow}
\HH^{k}(X, \partial X; \R) \rightarrow \ldots,
\end{equation}
$i^* \HH^{k-1}(X; \R)$ is a maximal isotropic subspace of
$\HH^{k-1}(\partial X; \R)$. 
Representing $\HH^{k-1}(\partial X; \R)$ by harmonic
forms, $\gamma \: i^* \HH^{k-1}(X; \R)$ is orthogonal to
$i^* \HH^{k-1}(X; \R)$.
By assumption, $\omega$ is orthogonal to
$i^* \HH^{k-1}(X; \R)$ and 
$\gamma \: i^* \HH^{k-1}(X; \R)$. Thus $\omega \: = \: 0$.
\end{proof}

Define $F^\prime_{\partial X} \in B(H_{\partial X})$ by
\begin{equation} \label{4.84}
F^\prime_{\partial X}(\omega) \: = \: 
\begin{cases}
\omega & \text{ if } \omega \in \overline{i^* \Ker(d)}, \\
- \: \omega & \text{ if } \omega \in \gamma \: \overline{i^* \Ker(d)}.
\end{cases}
\end{equation}

\begin{proposition} \label{4.85}
The triple $(H_{\partial X}, \gamma, F^\prime_{\partial X})$ 
represents
the same class in $\KKK_{2k-2}(C(\partial X); \C)$ as the triple
$(H_{\partial X}, \gamma, F)$ of Section \ref{Subsubsection 4.2.1}.
\end{proposition}
\begin{proof}
As $\HH^{k-1}(\partial X; \C)$ is finite-dimensional, 
$F^\prime_{\partial X} \: - \: F$ is compact.
\end{proof}

Proposition \ref{4.85} shows the K-cycle on $\partial X$
constructed from $X$, namely 
$(H_{\partial X}, \gamma, F^\prime_{\partial X})$, represents
the desired
K-homology class on $\partial X$.

\subsection{Even-dimensional quasiFuchsian manifolds}
\label{Subsubsection 4.5.2}

In this subsection we apply the formalism of Section
\ref{Subsubsection 4.5.1} to
describe an equivariant K-cycle on the limit set
of an even-dimensional 
quasiFuchsian manifold, in analogy with Section 
\ref{Subsubsection 4.4.3}.

We first consider the case of a Fuchsian manifold.
Let $\Gamma^\prime$ be a discrete torsion-free subgroup of 
$\Isom^+(H^{2k-1})$ whose limit set is $S^{2k-2}$.
There is a natural embedding $\Gamma^\prime \subset \Isom^+(H^{2k})$, 
with limit set $\Lambda^\prime \: = \: S^{2k-2} \subset S^{2k-1}$.
Applying Section \ref{Subsubsection 4.5.1} with $X$ being the
upper hemisphere $H^{2k-1} \subset S^{2k-1}$ gives
the K-cycle for
$\KKK_{2k-2}^{\Gamma^\prime}(C(S^{2k-2}); \C)$
of Section \ref{Subsubsection 4.2.4}.

A group $\Gamma \subset \Isom^+(H^{2k})$ 
that is quasiconformally related to $\Gamma^\prime$
is said to be a quasiFuchsian deformation of $\Gamma^\prime$.
Motivated by Section \ref{Subsubsection 4.4.3},
we can {\em define} a cycle for
$\KKK_{2k-2}^\Gamma(C(\Lambda); \C)$ by the pushforward under
$\phi \big|_{S^{2k-2}}$ of the K-cycle for
$\KKK_{2k-2}^{\Gamma^\prime}(C(S^{2k-2}); \C)$.
As in Section \ref{Subsubsection 4.4.3}, this is
independent of the choice of $\phi$. 
From Section \ref{Subsubsection 4.2.4},
the signature class for $S^{2k-2}$ is nontorsion in 
$\KKK_{2k-2}^{\Gamma^\prime}(C(S^{2k-2}); \C)$. As $\left( 
\phi \big|_{S^{2k-2}} \right)_*$ is an isomorphism, it follows that
the class in $\KKK_{2k-2}^\Gamma(C(\Lambda); \C)$ is also nontorsion.

\subsection{Even-dimensional convex-cocompact manifolds}
\label{Subsubsection 4.5.3}

In this subsection we apply the formalism of Section
\ref{Subsubsection 4.5.1} to
describe an equivariant K-cycle on the limit set
of a quasiconformal deformation of an even-dimensional convex-cocompact 
hyperbolic manifold whose convex core has totally geodesic boundary.

Let $\Gamma^\prime$ be a convex-cocompact subgroup of 
$\Isom^+(H^{2k})$ whose convex core has totally geodesic boundary.
Let $C$ be a connected component of $\partial \overline{M}$. Then the 
preimage $X$ of $C$ in $\Omega$ is 
a union $\bigcup_{i=1}^\infty B_i$ of round balls in
$S^{2k-1}$ with disjoint closures. Put $Y_i \: = \: \partial
\overline{B_i}$. Then the limit set $\Lambda^\prime$ is the closure of
$\bigcup_{i=1}^\infty Y_i$.

The Hilbert space of Section
\ref{Subsubsection 4.5.1} becomes
$H \: = \: \bigoplus_{i=1}^\infty L^{2}(Y_i; \Lambda^{k-1})$.
Define $\gamma_i \in B(L^{2}(Y_i; \Lambda^{k-1}))$ as in 
(\ref{4.5}).
Put $\gamma \: = \: \bigoplus_{i=1}^\infty \gamma_i$.
The operator $F$ of (\ref{4.84}) becomes a direct sum
$F \: = \: \bigoplus_{i=1}^\infty F_i$ where 
$F_i \in B(L^{2}(Y_i; \Lambda^{k-1}))$ is as in
(\ref{4.12}).
An element $a \in C(\Lambda^\prime)$ acts diagonally on $H$ as
multiplication by $a_i \: = \: a \big|_{Y_i}$ on 
$L^{2}(Y_i; \Lambda^{k-1})$.

\begin{proposition}
$(H, \gamma, F)$ is a
cycle for $\KKK^{\Gamma^\prime}_{2k-2}(C(\Lambda^\prime); \C)$.
\end{proposition}
\begin{proof}
Given $a \in C(\Lambda^\prime)$, we must show that $[F, a]$ is compact.
Extending $a$ to $a^\prime \in C(S^{2k-1})$ and approximating the latter
by smooth functions, we may assume that $a^\prime$ is smooth.

We know that for each $i$, $[F_i, a_i]$ is compact.  It suffices to show that
$\lim_{i \rightarrow \infty} \parallel [F_i, a_i] \parallel \: = \: 0$.
Fixing a round metric on $S^{2k-1}$, let $\overline{a}_i$ be the average
value of $a_i$ on $Y_i$. Then $[F_i, a_i] \: = \: [F_i, a_i \: - \:
\overline{a}_i]$ and 
$\lim_{i \rightarrow \infty} \parallel a_i \: - \: \overline{a}_i 
\parallel \: = \: 0$, from which the proposition follows.
\end{proof}

Now let $\Gamma$ be a quasiconformal deformation of $\Gamma^\prime$.
We can construct a cycle for
$\KKK_{2k-2}^\Gamma(C(\Lambda); \C)$ as the pushforward of the
preceding K-cycle by 
$\phi \big|_{\Lambda^\prime}$.
As in Section \ref{Subsubsection 4.4.3}, this is
independent of the choice of $\phi$. 

\subsection{The case of a Cantor set}
\label{Subsubsection 4.5.4}

In this subsection we specialize Section \ref{Subsubsection 4.5.3}
to the case
$k \: = \: 1$.

Let $\Gamma \subset \Isom^+(H^2)$ be a convex-cocompact subgroup.
If $M \: = \: H^2/\Gamma$ is noncompact then it has
a convex core with totally geodesic boundary, and
$\Lambda$ is a Cantor set.
Let $C$ be a connected component of $\Omega/\Gamma$. Then its
preimage $X$ in $\Omega$ is a countable
disjoint union of open intervals $(b_i, c_i)$ in $S^1$, and 
$\Lambda$ is the closure of the endpoints $\{b_i, c_i\}_{i=1}^\infty$. 
We have $H \: = \: l^2 \left(
\{b_i, c_i\}_{i=1}^\infty \right)$.
Define $\gamma \in B(H)$ by saying that for each
$\omega \in H$ and each $i$, $(\gamma \omega)(b_i) \: = \:
- \: \omega(b_i)$ and $(\gamma \omega)(c_i) \: = \: \omega(c_i)$. 
As $\Ker(d)$ consists of locally constant functions on $X$, we obtain
$(F \omega)(b_i) \: = \: \omega(c_i)$ and 
$(F \omega)(c_i) \: = \: \omega(b_i)$.

Taking a direct sum over the connected
components $C$ gives the K-cycle $(H, \gamma, F)$ considered in
\cite[Proposition 21, Section IV.3.$\epsilon$]{Connes (1994)}.
(The cited reference discusses $(H, F)$ as an ungraded K-cycle.)

\section{$p$-summability} \label{Subsection 4.6}

In this section we show the $p$-summability of a certain Fredholm module
$({\mathcal A}, H, F)$ for sufficiently large $p$.

With reference to Section \ref{Subsubsection 4.5.2},
let ${\mathcal A}$ be the restriction of $\phi^*
C^\infty(S^{2k-1})$ to $S^{2k-2}$,
a subalgebra of $C(S^{2k-2})$. Then we have an even Fredholm module
$({\mathcal A}, L^2(S^{2k-2}; \Lambda^{k-1}), F)$ in the sense of
\cite[Chapter IV, Definition 1]{Connes (1994)}.

\begin{proposition} \label{4.86}
For sufficiently large $p$, $({\mathcal A}, L^2(S^{2k-2}; \Lambda^{k-1}), F)$  
is $p$-summable in the sense of
\cite[Chapter IV, Definition 3]{Connes (1994)}.
\end{proposition}
\begin{proof}
We claim that for $p$ large,
$[F, a]$ is in the $p$-Schatten ideal
for all $a \in {\mathcal A}$.
Given $x, y \in S^{2k-2} \subset \R^{2k-1}$, 
let $|x-y|$ denote the chordal distance
between them.
From \cite{Janson-Wolff (1982)}, it suffices to show that
\begin{equation} \label{4.87}
\int_{S^{2k-2} \times S^{2k-2}}
\frac{|a(x) - a(y)|^p}{|x-y|^{4k-4}} \: dx \: dy \: < \: \infty.
\end{equation}
(The statement of \cite{Janson-Wolff (1982)} is for operators on 
$\R^{2k-2}$ instead of $S^{2k-2}$. We can go from one to the
other by stereographic projection, using the conformally-invariant
measure $\frac{dx \: dy}{|x-y|^{4k-4}}$.)
As $\phi$ is a quasiconformal homeomorphism, it lies in the H\"older
space $C^{0,\alpha}$ for some $\alpha \in (0,1)$. Then
there is a constant $C > 0$ such that $|a(x) - a(y)|^p \: \le \:
C \: |x-y|^{\alpha p}$ for all $x, y \in S^{2k-2}$.
The claim follows for $p \: > \: \frac{2k-2}{\alpha}$.
\end{proof}

With reference to Section \ref{Subsubsection 4.4.3},
let ${\mathcal A}$ be the restriction of $\phi^*
C^\infty(S^{2k})$ to $S^{2k-1}$,
a subalgebra of $C(S^{2k-1})$.
Let $E_{\pm}$ be the projection from
$L^2(S^{2k-1}; \Lambda^{k})$ to the $\pm 1$-eigenspace of
$\sign((-i)^k \: d*)$ acting on $\Image(d) \subset
L^2(S^{2k-1}; \Lambda^{k})$. Explicitly,
\begin{equation} \label{4.88}
E_{\pm} \: = \:
\frac12 \: \left( I \: \pm \frac{(-i)^k \: d*}{\triangle^{1/2}} \right)
\: \frac{dd^*}{\triangle}.
\end{equation}

For the motivation for the next proposition,
we refer to \cite[Section 7]{Connes (1985)}.

\begin{proposition} \label{4.89}
For sufficiently large $p$, $[E_\pm, a]$ is in the $p$-Schatten
ideal of operators on $L^2(S^{2k-1}; \Lambda^k)$ for all
$a \in {\mathcal A}$.
\end{proposition}
\begin{proof}
The proof is the same as that of Proposition \ref{4.86}.
\end{proof}

We note that Proposition \ref{4.89} refers to $L^2(S^{2k-1}; \Lambda^k)$,
whereas it is the $H^{-1/2}$-space $\Image(d) \subset
H^{-1/2}(S^{2k-1}; \Lambda^k)$ that is M\"obius invariant.  We can 
consider $L^2(S^{2k-1}; \Lambda^k)$ to be a dense subspace
of $H^{-1/2}(S^{2k-1}; \Lambda^k)$. The orthogonal projection
$E^\prime_\pm$ from $H^{-1/2}(S^{2k-1}; \Lambda^k)$ to the $\pm 1$-eigenspace
of $\sign((-i)^k \: d*)$ acting on $\Image(d) \subset
H^{-1/2}(S^{2k-1}; \Lambda^{k})$, i.e. to
$\Image \left( \frac{I \pm T}{2} \right)$, is again given by the formula in
(\ref{4.88}). Although we do not show the $p$-summability of the
ungraded Fredholm module $({\mathcal A}, H^{-1/2}(S^{2k-1}; \Lambda^{k}), 
E^\prime_+ - E^\prime_-)$, Proposition \ref{4.89} suffices for making
sense of the cyclic cocycles of \cite[Section 7]{Connes (1985)} in our case.

In the case $k = 1$ of Proposition \ref{4.89},
\cite[Section IV.3.$\gamma$, Proposition 14]{Connes (1994)} has the
stronger statement that
\begin{equation} \label{4.90}
\delta(\Gamma) \: = \: \inf \{p \: : \: [E_{\pm}, a] \text{ is in the
$p$-Schatten ideal for all } a \in {\mathcal A} \}.
\end{equation}
We do not know if a similar statement holds for all $k$. Using
\cite{Janson-Wolff (1982)}, it reduces
to a question about the Besov regularity of
$\phi \big|_{S^{2k-1}}$. The proof in
\cite[Section IV.3.$\gamma$, Proposition 14]{Connes (1994)} uses
facts about holomorphic functions that are special to the case $k = 1$.
One can ask the same question in the setup of Proposition 
\ref{4.86}.

Again in the case $k=1$, 
\cite[Section IV.3.$\gamma$, Theorem 17]{Connes (1994)} expresses
the Patterson-Sullivan measure on the limit set
in terms of the Dixmier trace.

\end{document}